\documentclass[a4paper, fleqn]{amsart}
\usepackage[english]{babel}
\usepackage{lmodern}
\usepackage[T1]{fontenc}
\usepackage[latin1]{inputenc}
\usepackage{array}
\usepackage{amsfonts,amssymb,amsmath}
\usepackage{mathrsfs}
\usepackage{esint}
\usepackage[all]{xy}
\usepackage{enumitem}
\usepackage{a4wide}
\usepackage{concrete}
\usepackage{eulervm}
\usepackage{tikz}
\usepackage[bottom]{footmisc}
\usepackage[colorlinks=true, urlcolor=blue, linkcolor=blue, citecolor=blue, pdftex]{hyperref}

\usetikzlibrary{arrows,chains,matrix,positioning,scopes}
\makeatletter
\tikzset{join/.code=\tikzset{after node path={%
\ifx\tikzchainprevious\pgfutil@empty\else(\tikzchainprevious)%
edge[every join]#1(\tikzchaincurrent)\fi}}}
\makeatother
\tikzset{>=stealth',every on chain/.append style={join},
         every join/.style={->}}
\tikzstyle{labeled}=[execute at begin node=$\scriptstyle,
   execute at end node=$]

\numberwithin{equation}{section}
\newtheoremstyle{mytheoremstyle} 
    {4mm}                    
    {4mm}                    
    {\itshape}                   
    {6mm}                           
    {\scshape}                   
    {.}                          
    {0.5em}                       
    {}  

\theoremstyle{mytheoremstyle}

\newtheorem{df}{Definition}[section]
\let\olddf\df
\renewcommand{\df}{\olddf\normalfont}
\newtheorem{thm}[df]{Theorem}
\newtheorem{prop}[df]{Proposition}

\newtheorem{lem}[df]{Lemma}

\newtheorem{ex}[df]{Example}
\let\oldex\ex
\renewcommand{\ex}{\oldex\normalfont}
\newtheorem{rk}[df]{Remark}
\let\oldrk\rk
\renewcommand{\rk}{\oldrk\normalfont}
\newtheorem*{pr}{Proof}
\let\oldpr\pr
\renewcommand{\pr}{\oldpr\normalfont}

\renewcommand{\epsilon}{\theta}
\newcommand{\C}{\mathbb{C}}
\newcommand{\R}{\mathbb{R}}
\newcommand{\Z}{\mathbb{Z}}
\newcommand{\N}{\mathbb{N}}

\newcommand{\A}{\mathcal{A}}
\newcommand{\G}{\mathcal{G}}
\renewcommand{\H}{\mathcal{H}}
\newcommand{\I}{\mathcal{I}}
\renewcommand{\L}{\mathcal{L}}
\renewcommand{\S}{\mathcal{S}}
\newcommand{\PS}{\mathcal{PS}}
\newcommand{\D}{\mathcal{D}}

\newcommand{\E}{\mathcal{E}}
\newcommand{\U}{\mathcal{U}}
\newcommand{\T}{\mathcal{T}}
\newcommand{\V}{\mathcal{V}}

\newcommand{\X}{\mathcal{X}}
\newcommand{\Y}{\mathcal{Y}}

\newcommand{\Tr}{\mathrm{Tr}}
\newcommand{\tr}{\mathrm{tr}}

\newcommand{\Ind}{\mathrm{Ind}}
\newcommand{\End}{\mathrm{End}}
\newcommand{\ch}{\mathrm{ch}}

\newcommand{\Vect}{\mathrm{Vect}}

\newcommand{\ad}{\mathrm{ad}}
\newcommand{\Ad}{\mathrm{Ad}}

\newcommand{\vol}{\mathrm{vol}}

\renewcommand{\dim}{\mathrm{dim}}

\renewcommand{\Re}{\mathrm{Re}}
\renewcommand{\Im}{\mathrm{Im}}
\newcommand{\Res}{\mathrm{Res}}
\renewcommand{\vol}{\mathrm{vol}}

\newcommand{\Td}{\mathrm{Td}}
\renewcommand{\ch}{\mathrm{ch}}
\renewcommand{\det}{\mathrm{det}}
\newcommand{\HP}{\mathrm{HP}}
\newcommand{\Pf}{\mathrm{Pf}}

\newcommand*{\barint}{\mathop{\ooalign{$\displaystyle{\int}$\cr$-$}}}

\newcommand{\alg}{\mathrm{alg}}
\renewcommand{\\}{\vspace{2mm}}
\renewcommand{\i}{\textup{\textbf{i}}}
\newcommand{\dd}{\textup{\textbf{d}}}

\newcommand{\DD}{\textup{\textbf{D}}}
\newcommand{\Hom}{\mathrm{Hom}}
\newcommand{\Diff}{\mathrm{Diff}}
\renewcommand{\tilde}{\widetilde}
\newcommand{\eps}{\varepsilon}
\newcommand{\ah}{\hat{a}}
\newcommand{\Th}{\widehat{T}}
\newcommand{\Hh}{\widehat{\H}}

\newcommand{\Yh}{\widehat{\Y}}
\newcommand{\Omh}{\widehat{\Omega}}

\newcommand{\cinf}{C^{\infty}}
\newcommand{\cinfc}{C^{\infty}_c}
\newcommand{\psio}{\overline{\psi}}
\newcommand{\Id}{\mathrm{Id}}

\title{An equivariant index theorem for hypoelliptic operators}
\author{Denis Perrot, Rudy Rodsphon}
\address{Universit\'e de Lyon \\ CNRS UMR 5208 \\ Universit\'e Lyon 1 \\ Institut Camille Jordan \\  43, Bd du 11 novembre 1918, 69622 Villeurbanne Cedex, France}
\email{perrot@math.univ-lyon1.fr}
\address{Vanderbilt University \\ Department of Mathematics \\ 1326 Stevenson Center \\ Nashville TN 37240, United States}
\email{rudy.rodsphon@vanderbilt.edu}

\begin{document}

\maketitle

\begin{abstract}
Let $M$ be a foliated manifold and $G$ a discrete group acting on $M$ by diffeomorphisms mapping leaves to leaves. Then $G$ naturally acts by automorphisms on the algebra of Heisenberg pseudodifferential operators on the foliation. Our main result is an index theorem for hypoelliptic-type operators which belong to the crossed product of the Heisenberg pseudodifferential operators with the group $G$. As a corollary, we get a solution to Connes-Moscovici's transverse problem in arbitrary codimensions, by exhibiting an explicit formula in terms of characteristic classes of equivariant vector bundles over $M$, for the Chern-Connes character associated to their hypoelliptic signature operator. 
\end{abstract}

\vskip 1cm

\setlength{\parindent}{0mm}
\textsc{Keywords.} Cyclic cohomology, K-theory, Index theory, Pseudodifferential operators \\

\textsc{MSC.}  19D55, 19K56, 58J42, 46L87

\setlength{\parindent}{0mm}

\section{Introduction}

This paper deals with the index theory of a certain class of equivariant hypoelliptic operators on foliated manifolds. Our motivation partly comes from \emph{transverse} index theory on foliations, and more specifically from the transverse index problem introduced by Connes and Moscovici in \cite{CM1995}, which attacks such questions under the utmost minimal hypotheses. In particular, one does not want to assume the foliation to be Riemannian nor the existence of any holonomy-invariant transverse measure. It is possible to reformulate that question as an equivariant index problem for (pseudo)groups of foliated diffeomorphisms acting on a certain foliated manifold without further assumption, i.e. we do not require the group action to preserve any geometric structure like a Riemannian metric or a conformal structure, unlike in classical index theory. A consequence of such generality is that elliptic operators are not a good object of study anymore and have to be replaced by a certain family of hypoelliptic operators, if one wants to get a relevant index theory. For that, one has to make extensive use of ideas from Noncommutative Geometry on the one hand, and of the Heisenberg pseudodifferential calculus on the other hand. One main obstacle is that well-established methods like the combination of heat kernel and Getzler rescaling methods do not yield convincing results even in codimension one. Purely topological methods from K-theory do not seem to apply well either (an explanation will be given below). To circumvent this issue, Connes and Moscovici developed the cyclic cohomology of Hopf algebras in \cite{CM1998}. Such an algebra plays the role of a group of symmetries and allows one to reorganize the calculations, leading to interesting results which are quite explicit in low codimensions. Nonetheless, their approach is rather indirect and in arbitrary codimensions, the full determination of the characteristic classes corresponding to their transverse index relies on a conjecture of Gelfand and Fuchs, which is still open. \\

The point of view developed in this paper completely avoids the use of Hopf algebras and gives in fact a more general index theorem. Moreover, it leads as a corollary to a direct identification of Connes-Moscovici's transverse index in terms of characteristic classes in equivariant cohomology, without further assumption. \\

We let $(M, V)$ be a (not necessarily compact) foliated manifold. To such a foliation is canonically associated the algebra $\Psi_{H,c}(M)$ of (compactly-supported) \emph{Heisenberg} pseudodifferential operators. The latter is a modification of the algebra of classical pseudodifferential operators, in which the vector fields tangent to the leaves of the foliation are of order 1, while the transverse vector fields are of order $\leq 2$. Let $G\subset \Diff(M)$ denote a discrete group of diffeomorphisms on $M$ mapping leaves to leaves. Then $G$ naturally acts on the algebra $\Psi_{H,c}(M)$ by automorphisms, and the algebraic crossed product $\Psi_{H,c}(M)\rtimes G$ is defined. We write any element $P \in \Psi_{H,c}(M) \rtimes G$ as a sum
\[ P = \sum_{g \in G} P_g \otimes U_g \]
where only finitely many coefficients $P_g\in \Psi_{H,c}(M)$ are non-zero. Notice that $P$ is naturally represented as a linear operator $\sum P_g \circ U_g : C^{\infty}_c(M) \longrightarrow \cinf_c(M)$, where $U_g$ acts as the shift operator $U_g(f)(x) = f(x \cdot g)$, for every $f\in \cinfc(M)$ and $x \in M$. Thus $P$ is far from being pseudodifferential in general, but belongs to the larger class of Fourier integral operators. At least at the algebraic level (see however the work of Savin and Sternin \cite{SS}), the index theory of such operators amounts to describing the $K$-theory/cyclic homology excision maps associated to the short exact sequence of algebras 
\[ 0 \longrightarrow \Psi^{-\infty}_c(M) \rtimes G \longrightarrow \Psi_{H,c}^{0}(M) \rtimes G \longrightarrow \S_{H,c}^{0}(M) \rtimes G \longrightarrow 0, \]
where $\Psi^{-\infty}_c(M)\subset \Psi_{H,c}^{0}(M)$ is the two-sided ideal of smoothing operators in the algebra of order $\leq 0$ Heisenberg pseudodifferential operators, and $\S_{H,c}^{0}(M) = \Psi_{H,c}^{0}(M)/\Psi^{-\infty}_c(M)$ is the quotient algebra of formal Heisenberg symbols. Let $\Tr_{[1]}$ be the trace on $\Psi_{c}^{-\infty}(M) \rtimes G$ obtained from the usual trace on $\Psi_{c}^{-\infty}(M)$ by localization at the unit of $G$ :
\[ \Tr_{[1]}\left(\sum_{g \in G} P_g \otimes U_g \right) = \Tr(P_1)\ . \]
We will mainly be interested in the image of $\Tr_{[1]}$ under the excision map in periodic cyclic cohomology $\partial: \HP^0(\Psi_{c}^{-\infty}(M) \rtimes G) \to \HP^1(\S_{H,c}(M) \rtimes G)$.

\begin{thm} \label{intro thm equivariant Radul cocycle} The boundary of the localized operator trace $\partial ([\Tr_{[1]}]) \in \HP^1(\S_{H,c}(M) \rtimes G)$ is represented by the \emph{equivariant Radul cocycle}
\[ \phi(P_g\otimes U_g , Q_h\otimes U_h) = \barint P_gU_g  [\ln \Delta_H^{1/4}, Q_hU_h] \quad \mbox{whenever} \ gh=1\ , \]
where $\Delta_H$ is the hypoelliptic sub-Laplacian (Example \ref{subelliptic laplacian}) associated to $M$, and the integral denotes the Connes-Moscovici residue over the algebra of formal symbols $\S^0_{H,c}(M)$.
\end{thm}

We refer to Section \ref{sindex} for a description of the Connes-Moscovici residue. The cocycle $\phi$ may be viewed as a \emph{local} formula in the sense that it only involves the formal Heisenberg symbol of pseudodifferential operators, and is given in terms of an integral over the \emph{Heisenberg cosphere bundle} $S^*_HM$ over $M$. A variant of this cocycle may also be obtained by taking the boundary of the localized operator trace under the extension
\[ 0 \longrightarrow \Psi_{H,c}^{-1}(M) \rtimes G \longrightarrow \Psi_{H,c}^{0}(M) \rtimes G \longrightarrow \cinfc(S^{*}_H M) \rtimes G \longrightarrow 0\ , \]
leading to a cocycle $\partial([\tau])\in \HP^1 (\cinfc(S^{*}_H M) \rtimes G)$ over the algebra of Heisenberg principal symbols, which may easily be related to $\phi$. The main result of this paper (Theorem \ref{ttodd}) is a geometric realization of this cocycle :

\begin{thm}  \label{intro thm equivariant index}
Let $M$ be a foliated manifold and $G$ be a discrete group of foliated diffeomorphisms. Let $EG$ be the universal bundle over the classifying space $BG$ of $G$. Let 
\[ 0 \longrightarrow \Psi^{-1}_{H,c}(M)\rtimes G \longrightarrow \Psi^0_{H,c}(M)\rtimes G \longrightarrow C^{\infty}_c (S^*_HM) \rtimes G \longrightarrow 0 \] be the equivariant Heisenberg pseudodifferential extension. Then, the image of the localized trace at the unit $\partial([\tau])\in \HP^1 (\cinfc(S^{*}_H M) \rtimes G)$ by excision is given by
\begin{equation*}
\partial([\tau]) = \Phi(\Td(TM\otimes\C)) 
\end{equation*}
where $\Phi: H^{\mathrm{ev}}(EG\times_{G} S^*_HM) \to HP^1(C^{\infty}_c(S^*_HM)\rtimes G)$ is Connes' characteristic map from equivariant cohomology to cyclic cohomology, and $\Td(TM\otimes\C) $ is the equivariant Todd class of the complexified tangent bundle of $M$. 
\end{thm}

We stress once again that this result holds for \emph{any} group of foliated diffeomorphisms $G$. Let us now explain the novelties of our approach and why it allows such flexibility. \\

The machinery developed in this paper is based upon an adaptation of the formalism developed in \cite{Per2012} and \cite{Per2014} to the Heisenberg calculus. In the non-foliated case, this formalism can be viewed as a cohomological counterpart to the K-homological approach to index theory developed by Kasparov (cf. for example \cite{Kas}). The main point of his approach is, after a choice of an almost complex (hence $\text{Spin}^c$) structure on the cotangent bundle, the construction of a fundamental K-homology class using the Dolbeault-Dirac operator on $T^*M$. As long as one deals with isometric actions, Kasparov's approach extends to the equivariant case without much complications, and the Chern character of the Dolbeault-Dirac element can be calculated by standard methods. However, it is not possible to work with arbitrary actions : if one wants to get an invariant Dolbeault-Dirac operator, the lifted action of the group on $T^*M$ must at least preserve the almost complex structure chosen. \\

A fundamental idea of \cite{Per2012} used here is the replacement of the Dolbeault-Dirac operator by a more canonical object that we call a \emph{generalized Dirac operator} (cf. Section \ref{Dirac operators}). Heuristically, it is constructed from a Dirac operator on $T^*M$, but uses directly the canonical symplectic structure of $T^*M$ instead of its almost complex structure, and a $\text{Spin}(n,n)$-structure instead of a $\text{Spin}^c$-structure (see Section \ref{spinors} for details). A significant benefit is that now, the lifted action to $T^*M$ of any group of diffeomorphisms on $M$ preserves these structures. This can be seen at the level of structure groups : the general linear group $GL_n(\R)$ embeds in both the symplectic group $Sp(2n)$ and the indefinite orthogonal group $O(n,n)$. The next step is to carry the action of these Dirac operators at a quantized level, that is on symbols rather than on functions on the cotangent bundle, which gives rise to generalized Dirac operators. In a sense, an interesting feature of this formalism is to develop a purely \emph{algebraic} heat kernel method with a hyperbolic laplacian $\partial_x \partial_p$, where $(x,p)$ is a coordinate system on $T^*M$. This would not be feasible if one was not working at the symbolic level. Let us mention as well that this strategy leads to a simpler and more direct approach to the Nest-Tsygan deformation quantization approach to index theory \cite{NT}. This relationship will be developed in more detail elsewhere.  \\
 
Finally, one shows that a periodic cyclic cohomology class, which computes the index, can be associated to such an operator. This is done by constructing an algebraic JLO formula, which is an adaptation of the usual ingredients of the JLO cocycle \cite{JLO} to this context. Hence, this class may be viewed as a formal analogue of the Chern character in K-homology. The fact that generalized Dirac operators are \emph{diffeomorphism invariant} finally explains why we can extend these constructions to the equivariant setting without any restriction on the group action, following \cite{Per2014}. All these ideas can be set up without major problems to foliations and Heisenberg symbols. Some new subtleties appear in the present paper though : since we do not suppose that $EG$ be a manifold here, we have to deal with certain technical issues. One of them is that we have to work with with non-metric connections. \\

Back to foliations, Theorem \ref{intro thm equivariant index} is proved by computing the JLO cocycle for two different Dirac operators : the first one gives the Radul cocycle of the pseudodifferential extension, while the second one gives the equivariant Todd class. Usual homotopy arguments from the original JLO formula apply verbatim there, and allows us to conclude the proof. Contrary to what happens in the ordinary JLO setup, the calculations are completely algebraic, because our JLO cocycles are built from a Wodzicki-type residue and therefore directly extract terms up to a finite order in formal symbols. \\

We then apply Theorem \ref{intro thm equivariant index} to the Connes-Moscovici index problem for \emph{transversally} hypoelliptic operators on foliations. After reduction to a complete transversal $W$, the holonomy groupoid of a given foliation is Morita equivalent to an \'etale groupoid $W\rtimes G$ where $G\subset \Diff(W)$ is a discrete (pseudo)group of diffeomorphisms. The main object of study is therefore the crossed product algebra $\cinfc(W) \rtimes G$ where, for notational simplicity, we treat $G$ as a group. Since $G$ is not supposed to preserve any geometric structure on $W$, $G$-invariant elliptic differential operators do not exist even at the leading symbol level. The idea of Connes in \cite{C83} is to pass by a Thom isomorphism to the bundle of Riemannian metrics over $W$. This fibration is in particular a foliation $M$, \emph{to which Theorem \ref{intro thm equivariant index} will be applied}. The action of $G$ on $W$ lifts to $M$, mapping leaves to leaves. Connes and Moscovici then construct a hypoelliptic signature operator on $M$ which is almost invariant under the $G$-action, in the sense that its Heisenberg leading symbol is $G$-invariant. This yields a regular spectral triple over the algebra $\cinfc(M)\rtimes G$ whose Chern-Connes character, which represents the index of the Connes-Moscovici signature operator, may be computed by means of a Wodzicki-type residue formula (\cite{CM1995}). However, this does not directly provide a characteristic class formula, since the actual calculations give thousands of terms\footnote{Whereas it is directly computable when applied in the setup of classical index theory, by Getzler's rescaling method for example.} already for very low-dimensional manifolds $W$. To overcome this difficulty in higher dimensions, Connes and Moscovici used Hopf algebras and their cyclic cohomology in \cite{CM1998} to prove, under the mild assumption that the lifted action of $G$ to $M$ has no fixed points, that the Chern-Connes character of the hypoelliptic signature operator can be realized in Gelfand-Fuchs cohomology. They predict that the corresponding class is a universal polynomial in the Pontryagin classes. Explicit calculations are made in \cite{CM1998} in dimension 1, giving (twice) the transverse fundamental class of \cite{C83}. In dimension 2, the authors show that the coefficient of the first Pontryagin class does not vanish. \\

Our Theorem \ref{intro thm equivariant index} allows us to shortcut the calculation with Hopf algebras and gives a direct answer to the problem of computing the Chern-Connes character of the hypoelliptic signature operator in terms of equivariant characteristic classes, for manifolds $W$ of arbitrary dimension. This answers to Connes-Moscovici's affrimation positively. 

\begin{thm} \label{intro thm CM signature}
Let $G$ be a discrete group of orientation-preserving diffeomorphisms on a manifold $W$. Let $M$ be the bundle of Riemannian metrics over $W$. If the lifted action of $G$ has no fixed points on $M$, then the Chern-Connes character of the Connes-Moscovici hypoelliptic signature operator is   
\[ \pi_*\circ\Phi(L'(M))\ \in HP^1(C^{\infty}_c(M)\rtimes G)\ , \]
where $L'(M)$ is a modified $L$-genus, $\Phi: H^{\mathrm{ev}}(EG\times_G S^*_HM) \to HP^1(C^{\infty}_c(S^*_HM)\rtimes G)$ is Connes' characteristic map, and $\pi_*: HP^1(C^{\infty}_c(S^*_HM)\rtimes G) \to HP^1(C^{\infty}_c(M)\rtimes G)$ is induced by the canonical projection $\pi : S^*_HM\to M$.
\end{thm}

\emph{Organization of the paper.} Section \ref{sindex} is a self-contained introduction to the equivariant Heisenberg pseudodifferential calculus on foliated manifolds, the pseudodifferential extension, and the Connes-Moscovici residue. Then Theorem \ref{intro thm equivariant Radul cocycle} is proved. \\

In Sections \ref{sbimod}, \ref{scano} and \ref{sdirac} we adapt the formalism of \cite{Per2012} to the Heisenberg calculus. We introduce various spaces of operators acting on formal Heisenberg symbols, and the universal Dirac operators which will be used in our algebraic JLO formula. \\

Section \ref{sequiv} introduces the required objects to carry this formalism to the equivariant setting. In particular, we recall the point of view we need to construct Connes' characteristic map from the equivariant cohomology $H^{\bullet}(EG \times_G M)$ to the periodic cyclic cohomology of the crossed product $\HP^{\bullet}(\cinf(M) \rtimes G)$. From the technical side we use the $X$-complex of Cuntz and Quillen \cite{CQ95}. \\

Section \ref{sJLO} finally gives the algebraic JLO formula on the algebra of formal (equivariant) Heisenberg symbols, leading to Theorem \ref{intro thm equivariant index}. This is again an adaptation of the formalism developed in \cite{Per2014} for the non-Heisenberg case.\\ 

Section \ref{CM signature} shows how to deduce Theorem \ref{intro thm CM signature} from Theorem \ref{intro thm equivariant index}. \\

\section*{General notations}

\subsection*{Differential geometry} Throughout the paper, $M$ will be a (not necessarily compact) foliated manifold without boundary of dimension $n=p+q$, $q$ being the codimension of the foliation, $G$ will denote a discrete group of foliated diffeomorphisms. The only exception is Section \ref{sequiv}, where $M$ will be any manifold without boundary, $G$ any discrete group of orientation-preserving diffeomorphisms and $n$ does not necessarily denote the dimension of $M$. The bundles $TM$, $T^*M$ and $S_H^* M$ will respectively denote the tangent, cotangent and Heisenberg cosphere bundles of $M$. The spaces $BG$ and $EG$ are respectively the classifying space of $G$ and the universal $G$-bundle over it, as usual. We shall denote $d_H$ the de Rham differential on $EG$ (seen as a simplicial manifold). 

\subsection*{Spinors} We will denote $E$ the complexified exterior powers $\Lambda^{\bullet} T^*M \otimes \C$ of $T^*M$, which is bundle of spinors over $M$ for the Clifford algebra of $TM \oplus T^*M$, the quadratic form on $TM \oplus T^*M$ considered is the duality bracket (and has signature $(n,n)$). Locally, a set of generators of $\End(E)$ is given by $\psi^i, \psio_i, i=1, \ldots, n$, which respectively denote the exterior product with $dx^i$ and the interior product with $\partial_{x^i}$. 

\subsection*{Pseudodifferential calculus} $\Psi^{-\infty}(M)$, $\Psi_H(M)$ and $\S_H(M)$ will respectively denote the algebra of smoothing operators, classical Heisenberg pseudodifferential operators and classical formal symbols on $M$. We denote $\S(M)=\L(M)[[\eps]]$ (formal series), where $\L(M)$ is the bimodule of formal symbols on $M$, and $\D(M)$ a subalgebra of $\S(M)$ that contains the operators of interest (cf. Section \ref{sbimod} for precise definitions). An added subscript $c$ stands for compactly supported version of these spaces. Versions of these algebras acting on sections of vector bundles will also be used.  

\subsection*{Algebras} For any associative algebra $\A$ over $\C$, let $T\A$ and $J\A$ denote respectively the tensor algebra of $\A$ and the kernel of the multiplication map $T\A \to \A$. The $J$-adic completion $\varprojlim_m T\A/(J\A)^m $ of $T\A$ will be denoted $\Th \A$. 

\subsection*{Cyclic (co)homology} Let $\A$ be any associative algebra over $\C$. We shall make extensive use of the (completed) space of non-commutative differential forms $\Omh\Th\A=\prod_{n=0}^{\infty}\Omega^n\Th\A$ endowed with the total differential $(b+B)$, which computes the periodic cyclic homology of $\A$. Another complex that we use to get periodic cyclic cocycles is the $X$-complex $X(\A)$ of Cuntz-Quillen. \\

\section{Equivariant local index formula}\label{sindex}

\subsection{Heisenberg pseudodifferential calculus on foliations}

Let $M$ be a foliated manifold of dimension $n$, and let $V$ be the integrable sub-bundle of the tangent bundle $TM$ of $M$ which defines the foliation. Denote $v$ the dimension of the leaves and $h = n-v$ their codimension.  \\

The fundamental idea of the Heisenberg calculus is that \emph{longitudinal} vector fields (with respect to to the foliation) have \emph{order 1}, whereas \emph{transverse} vector fields have \emph{order $\leq$ 2}. We shall now describe the symbolic calculus in more detail, following Connes and Moscovici \cite{CM1995}.  \\

Let $(x_1, \ldots , x_n)$ be a foliated local coordinate system of $M$, i.e, the vector fields $\frac{\partial}{\partial x_1}, \ldots , \frac{\partial}{\partial x_v}$ locally span $V$, so that $\frac{\partial}{\partial x_{v+1}}, \ldots , \frac{\partial}{\partial x_n}$ are transverse to the leaves of the foliation. Then, we set 
\begin{align*}
& \vert p \vert' = (p_1^4 + \ldots + p_v^4 + p_{v+1}^2 + \ldots + p_{n}^2)^{1/4} \\ & \langle \alpha \rangle = \alpha_1 + \ldots + \alpha_v + 2\alpha_{v+1} + \ldots 2\alpha_n 
\end{align*}
for every $p \in \R^n$, $\alpha \in \N^n$.

\begin{df} A smooth function $a(x,p) \in C^\infty(\R^n_x \times \R^n_p)$ is a \emph{Heisenberg symbol of order $m \in \R$} if over any compact subset $K\subset\R^n_x$ and for every multi-index $\alpha, \beta$, one has the following estimate
\begin{equation*} \vert \partial^\beta_x \partial^\alpha_p a(x,p) \vert \leq C_{K,\alpha,\beta}(1+\vert p \vert')^{m - \langle \alpha \rangle} \end{equation*}
\end{df}

We shall focus on the smaller class of \emph{classical Heisenberg symbols}. For this, we first define the \emph{Heisenberg dilations}
\begin{equation*} 
\lambda \cdot (p_1, \ldots, p_v, p_{v+1}, \ldots, p_n) = (\lambda p_1, \ldots, \lambda p_v, \lambda^2 p_{v+1}, \ldots, \lambda^2 p_n) 
\end{equation*}
for any non-zero $\lambda \in \R_+$ and non-zero $p \in \R^n$. 

Then, a Heisenberg pseudodifferential symbol $\sigma$ of order $m$ is called \emph{classical} if it has an asymptotic expansion when $\vert p \vert' \to \infty$
\begin{equation} \label{classical hpdo}
a(x,p) \sim \sum_{j \geq 0} a_{m-j}(x,p) 
\end{equation} 
where $a_{m-j}(x,p)$ are \emph{Heisenberg homogeneous} functions, that is, for any non zero $\lambda \in \R$,
\begin{equation*} a_{m-j}(x,\lambda \cdot p) = \lambda^{m-j} a_{m-j}(x, p) \end{equation*} 
The \emph{Heisenberg principal symbol} is the symbol of highest degree in the expansion (\ref{classical hpdo}). \\ 

To such a symbol $a$ of order $m$, one associates its left-quantization as the linear map:
\begin{equation*} 
P : C^\infty_c(\R^n) \to C^\infty(\R^n), \qquad Pf(x) = \frac{1}{(2\pi)^n} \int_{\R^n} e^{\i x \cdot p} a(x,p) \hat{f}(p) dp 
\end{equation*}
where $\hat{f}$ denotes the Fourier transform of the function $f$. We shall say that $P$ is a \emph{classical Heisenberg pseudodifferential operator of order $m$}. If $P$ is properly supported, then it actually defines a linear map $C^\infty_c(\R^n) \to C^\infty_c(\R^n)$. We denote by $\Psi_H^m(\R^n)$ the vector space of such properly-supported operators and by $\Psi_{H,c}^m(\R^n)$ its subspace of \emph{compactly} supported operators. Since properly-supported operators can be composed, the unions of all-orders operators
\begin{equation*} \Psi_H(\R^n) = \bigcup_{m \in \R} \Psi_H^m(\R^n)\ ,\qquad  \Psi_{H,c}(\R^n) = \bigcup_{m \in \R} \Psi_{H,c}^m(\R^n)
\end{equation*}
are associative algebras over $\C$. The ideals of \emph{regularizing operators}
\begin{equation*} \Psi^{-\infty}(\R^n) = \bigcap_{m \in \R} \Psi_H^m(\R^n) \ ,\qquad         \Psi^{-\infty}_c(\R^n) = \bigcap_{m \in \R} \Psi_{H,c}^m(\R^n) 
\end{equation*}
correspond respectively to the algebras of operators with properly- and compactly-supported smooth Schwartz kernel.  \\

If $P_1,P_2 \in \Psi_H(\R^n)$ are Heisenberg pseudodifferential operators of symbols $a_1$ and $a_2$, $P_1P_2$ is a Heisenberg pseudodifferential operator whose symbol $a$ is given by the \emph{star-product} of symbols :
\begin{equation} \label{star product} 
a(x, p) = a_1 \star a_2 (x,p) \sim \sum_{\vert \alpha \vert \geq 0} \frac{(-\i)^{\vert \alpha \vert}}{\alpha!} \partial^\alpha_p a_1(x, p) \partial^\alpha_x a_2(x, p) 
\end{equation} 
Notice that the order of each symbol in the sum is decreasing while $\vert \alpha \vert$ is increasing. \\

We define the \emph{algebra of Heisenberg formal classical symbols} $\S_H(\R^n)$ and its compactly-supported subalgebra $\S_{H,c}(\R^n)$ as quotients
\begin{equation*}
\S_H(\R^n) = \Psi_H(\R^n)/\Psi^{-\infty}(\R^n) \ ,\qquad \S_{H,c}(\R^n) = \Psi_{H,c}(\R^n)/\Psi_c^{-\infty}(\R^n)
\end{equation*}
Their elements are formal sums given in (\ref{classical hpdo}), and the product is the star product (\ref{star product}). \\ 

A Heisenberg symbol is said \emph{Heisenberg elliptic} if it is invertible in $\S_H(\R^n)$. This is equivalent to say that its \emph{Heisenberg principal symbol} is invertible on $\R^n_x \times \R^n_p \smallsetminus \{0\}$. The corresponding pseudodifferential operators are in general not elliptic, but only hypoelliptic. 

\begin{ex} \label{subelliptic laplacian} The hypoelliptic sub-Laplacian is the differential operator  
\begin{equation*}
\Delta_H = \partial_{x_1}^4 +  \ldots + \partial_{x_v}^4 - (\partial_{x_{v+1}}^2 + \ldots + \partial_{x_n}^2)
\end{equation*} 
It has Heisenberg principal symbol $a(x,p) = \vert p \vert'^4$, and is therefore Heisenberg elliptic. However, its usual principal symbol, as an ordinary differential operator, is $(x,p) \mapsto \sum_{i=1}^v p_i^4 $, so $\Delta_H$ is clearly not elliptic.  \\
\end{ex}

Heisenberg pseudodifferential operators are compatible with foliated coordinate changes. Therefore, the Heisenberg calculus can be defined globally on foliations by using a partition of unity. Then, for a foliated manifold $M$, we denote by $\Psi_H(M)$ the algebra of properly-supported Heisenberg pseudodifferential operators on $M$, and by $\Psi_{H,c}(M)$ its subalgebra of compactly-supported operators. \\

For a ($\Z_2$-graded) complex vector bundle $E$ over $M$, one defines in the same way the algebra of Heisenberg pseudodifferential operators $\Psi_H(M,E)$ acting on the smooth compactly-supported sections $C^{\infty}_c(M,E)$ of $E$. One always has an exact sequence 
\[ 0 \to \Psi^{-\infty}(M,E) \to \Psi_H(M,E) \to \S_H(M,E) \to 0 \]
and similarly for the algebra of compactly-supported operators $\Psi_{H,c}(M,E)$. Notice that for $a \in \S_H(M,E)$, $(x,p) \in T^{*}_xM$, we have $a(x,p) \in \End(E_x)$. Let $\PS_H(M,E) \subset \S_H(M,E)$ denote the subalgebra of polynomial Heisenberg symbols (with respect to the cotangent coordinate $p$). The latter is isomorphic to the algebra of differential operators, endowed with the Heisenberg degree. \\

\subsection{Spinors and commutation relations} \label{spinors}
The vector bundle of interest in this paper will be the exterior algebra $E=\Lambda^{\bullet} (T^*M\otimes\C)$ of the complexified cotangent bundle. In a foliated coordinate system $(x^{1}, \ldots, x^{n})$ over an open subset $U \subset M$, a local basis of the sections of $E$ is given by $ 1,dx^{i_1}, dx^{i_1} \wedge dx^{i_2}, \ldots, dx^{i_1} \wedge \ldots \wedge dx^{i_n}$, $1 \leq i_1< \ldots <i_n \leq n $. Moreover, the endomorphisms $\End(E_x)$ of the fibre $E_x$ are generated by   
\begin{equation}
\psi^{i} = dx^{i} \wedge \textbf{.}, \quad \overline{\psi}_i = \iota(\partial_{x^{i}})
\end{equation}
for $i = 1, \ldots, n$, where $\iota$ stands for the interior product with a vector field. We have the following \emph{anti-commutation relations rules} :  
\begin{equation}
[\psi^{i}, \overline{\psi}_j] = \delta_i^{j}, \quad [\psi^{i}, \psi^{j}] = [\overline{\psi}_i, \overline{\psi}_j] = 0
\end{equation}
Then, $\End(E_x)$ is the Clifford algebra of $T_x M \oplus T^{*}_x M$, the quadratic form is the duality bracket, which has signature $(n,n)$. In other words, $E$ is a bundle of spinors for the Clifford algebra bundle $\text{Cliff}(TM \oplus T^{*}M)$ over $M$. Let us also recall the commutation relations of symbols : in the coordinate system $(x,p)$ over $T^{*}U$, we have : 
\begin{equation}
[x^{i},p_j] = -\i \delta^i_{j}, \quad [x^{i},x^{j}] = [p_i, p_j] = 0
\end{equation} 
Thus, every Heisenberg symbol $a \in \S_H(M,E)$ is locally over $M$ a formal series  
\[ a(x,p,\psi,\overline{\psi}) = \sum_{j=0} a_{m-j}(x,p,\psi,\overline{\psi}) \]
where the functions $a_{m-j}$ are Heisenberg-homogeneous in $p$ and polynomial in the variables $\psi$ and $\overline{\psi}$. 

\subsection{Wodzicki residue on $\Psi_H(M)$}
Let $M$ be a foliated manifold. By choosing a Riemannian metric on $M$, one can construct a hypoelliptic sub-Laplacian $\Delta$ as in the flat example \ref{subelliptic laplacian}. Throughout this article, we define the complex powers $\Delta^{-z}$ as properly-supported Heisenberg pseudodifferential operators, using a parametrix $(\lambda - \Delta)^{-1}$ and an appropriate Cauchy integral 
\[ \Delta^{-z} = \dfrac{1}{2 \pi \i} \int \lambda^{-z} (\lambda - \Delta)^{-1} \, d\lambda \]
where the contour is a vertical line pointing downwards. $\Delta^{-z}$ is uniquely defined modulo addition of a smoothing operator.\\

\begin{thm} \label{CMzeta} \emph{(Connes-Moscovici, \cite{CM1995})} Let $M$ be a foliated manifold of dimension $n$, $v$ be the dimensions of the leaves, $h$ their codimension, and $P \in \Psi^m_{H,c}(M)$ be a compactly-supported Heisenberg pseudodifferential operator of order $m \in \R$. Then, for any sub-elliptic sub-laplacian $\Delta$, the zeta function 
\begin{equation*} \zeta_P(z) = \Tr(P \Delta^{-z/4}) \end{equation*}
is holomorphic on the half-plane $\Re(z) > m+v+2h$, and extends to a meromorphic function over the whole complex plane, with at most simple poles in the set 
\begin{equation*} \{ m+v+2h, m+v+2h - 1, \, \ldots \} \end{equation*}
\end{thm}

The meromorphic extension of the zeta function given by this theorem allows the construction of a Wodzicki-Guillemin trace on $\S_{H,c}(M) = \Psi_{H,c}(M)/\Psi_c^{-\infty}(M)$.

\begin{thm} \emph{(Connes-Moscovici, \cite{CM1995})} \label{CM residue} The Wodzicki residue functional 
\begin{equation*}  
\barint : \Psi_{H,c}(M) \longrightarrow \C , \quad P \longmapsto \Res_{z=0} \Tr(P \Delta^{-z/4}) 
\end{equation*}
is a trace vanishing on the ideal of smoothing operators, hence descends to a trace on $\S_{H,c}(M)$. Moreover we have the following formula, only depending on the formal symbol $a$ of $P$ up to finite order: 
\begin{equation}  \label{wodzicki residue symbol 0}
\barint P = \frac{1}{(2\pi)^n} \int_{S_H^*M} \iota_L \left(a_{-(v+2h)}(x,p) \, \frac{ \omega^n}{n!}\right) 
\end{equation}
\end{thm}

Here, $S_H^*M$ is the Heisenberg cosphere bundle, that is, the sub-bundle 
\begin{equation*} 
S_H^*M = \{(x,p) \in T^*M \, ; \, \vert p \vert' = 1 \} 
\end{equation*} 
of the cotangent bundle $T^*M$, $\omega$ denotes the standard symplectic form on $T^*M$, $\iota$ stands for the interior product, and $L$ is the generator of the Heisenberg dilations given by the formula
\[ L = \sum_{i=1}^v p_i \partial_{p_i} + 2\sum_{i=v+1}^n p_i \partial_{p_i} \]
Notice that the residue always vanishes on operators with non-integral order $m$. \\

All these results still hold for Heisenberg pseudodifferential operators acting on sections of a ($\Z_2$-graded) vector bundle $E$ over $M$. In this case, the symbol $a_{-(v+2n)}(x,p)$ above is at each point $(x, p)$ an endomorphism acting on the fibre $E_x$, and (\ref{wodzicki residue symbol 0}) becomes :
\begin{equation*} 
\barint P = \frac{1}{(2\pi)^n} \int_{S_H^*M} \iota_L \left(\tr _s(a_{-(v+2n)}(x,p)) \, \frac{ \omega^n}{n!}\right) 
\end{equation*}
where $\tr_s$ denotes the fibrewise (graded) trace of endomorphisms. \\

\subsection{Excision and equivariant residue index formula}

Let $M$ be a foliated manifold. Consider a discrete subgroup $G\subset\Diff(M)$ of diffeomorphisms mapping leaves to leaves. By convention we suppose that $G$ acts from the right, so, for any $g\in G$ the induced linear action $U_g$ on the space of functions $\cinf(M)$ reads
\[ (U_g f)(x) = f(x\cdot g)\ ,\quad \forall f\in\cinf(M)\ ,\ x\in M \]
Recall that the \emph{algebraic crossed product} $\Psi_{H,c}(M) \rtimes G$ is the universal algebra generated by Heisenberg pseudodifferential operators and group elements, that is,
\[ \Psi_{H,c}(M) \rtimes G = \left\{ \sum_{g \in G} P_g \otimes U_g \, ; \, P_g \in \Psi_{H,c}(M) \right\} \]
and the sum only contains a finite number of non-zero coefficients $P_g$. The multiplication is given by the rule 
\[ (P\otimes U_g) \cdot (Q\otimes U_h) = P (U_g Q U_{g^{-1}})\otimes U_{gh}  \]
To this effect, remark that $U_g Q U_{g^{-1}}$ is still a classical Heisenberg pseudodifferential operator, so that the product with $P$ makes sense. Note also that in general, the representation of $\Psi_{H,c}(M) \rtimes G$ as linear operators on $\cinf(M)$ does not yield pseudodifferential operators.  \\

Then, one has an \emph{extension}
\begin{equation}  \label{extension}
0 \to \Psi_c^{-\infty}(M) \rtimes G \to \Psi^0_{H,c}(M) \rtimes G \to \S^0_{H,c}(M) \rtimes G \to 0 
\end{equation}
where $\Psi^0_{H,c}(M)$ denotes the subalgebra of operators of \emph{integral} order $m\leq 0$. The usual trace on the algebra of regularizing operators $\Psi_c^{-\infty}(M)$, given by 
\[ \Tr(K) = \int_M k(x,x) d\vol(x) \]
where $k$ stands for the Schwartz kernel of $K$, extends to a trace $\Tr_{[1]}$ on $\Psi_c^{-\infty}(M) \rtimes G$ by \emph{localization at the unit} of $G$
\[ \Tr_{[1]}\left(\sum_{g \in G} K_g\otimes U_g \right) = \Tr(K_1) \]
That $\Tr_{[1]}$ still remains a trace only comes from the invariance of the ordinary operator trace $\Tr$ under conjugation by $U_g$. \\

In the same way, the Wodzicki residue on $\S_{H,c}(M)$ extends to a trace on $\S_{H,c}(M) \rtimes G$ by localization at the unit.   \\

The pseudodifferential extension gives rise to the following commutative diagram involving \emph{algebraic K-theory} and \emph{periodic cyclic homology} (\cite{Nis1997})

\begin{equation} \label{Index diagram}
\vcenter{
\xymatrix{
    K_1^{\alg}(\S^0_{H,c}(M) \rtimes G) \ar[r]^{\Ind} \ar[d]^{\ch_1} & K_0^{\alg}(\Psi_c^{-\infty}(M) \rtimes G) \ar[d]^{\ch_0} \\
    \HP_1(\S^0_{H,c}(M) \rtimes G) \ar[r]^{\partial} & \HP_0(\Psi_c^{-\infty}(M) \rtimes G)}
}
\end{equation}
The vertical arrows are respectively the odd and even Chern character. \\
 
Denote again $\partial :  \HP^0(\Psi_c^{-\infty}(M) \rtimes G) \rightarrow \HP^1(\S_{H,c}^0(M) \rtimes G)$ the induced excision map in cohomology. We shall now compute $\partial ([\Tr_{[1]}])$. To do this, we lift $\Tr_{[1]}$ on $\Psi_c^{-\infty}(M) \rtimes G$ to a linear map on $\Psi_{H,c}^0(M) \rtimes G$ using a zeta function renormalization 
\[ \Tr_{[1]}'\left(\sum_{g \in G} P_g\otimes U_g \right) = \Pf_{z=0} \Tr \left(P_1  \cdot \Delta^{-z/4}\right)  \]
where $\Delta$ is a sub-elliptic Laplacian, and $\Pf_{z=0}$ is the constant term in the Laurent series expansion of the zeta-function at $z=0$. Then,  $\partial ([\Tr_{[1]}])$ is represented in $\HP^1(\S_H^0(M) \rtimes G)$ by the cyclic 1-cocycle 
\begin{equation}
\phi(P\otimes U_g,Q\otimes U_h) = \Tr_{[1]}'([PU_g,QU_h])
\end{equation}
for all $P,Q\in \Psi^0_{H,c}(M)$ and $g,h\in G$. This expression makes sense as a cocycle over $\S_H^0(M) \rtimes G$ because it vanishes whenever $P$ or $Q$ belongs to the smooting ideal $\Psi^{-\infty}_{H,c}(M)$. Then, because the trace is localized at units, one finds that $\phi(P\otimes U_g,Q\otimes U_h)=0$ if $gh\neq 1$, and
\begin{equation*}  \phi(P\otimes U_g, Q\otimes U_{g^{-1}}) = \Tr_{[1]}'([P U_g, Q U_{g^{-1}}]) = \Pf_{z=0} \Tr \left([P U_g, Q U_{g^{-1}}]  \cdot \Delta^{-z/4}\right) \end{equation*}
otherwise. The formula can be made a little more explicit if we gather accurately the relevant terms. This is the aim of the following proposition. 

\begin{prop} The cyclic 1-cocycle $\phi$ given above is given in terms of the Connes-Moscovici residue :
\begin{equation*} \label{equivariant Radul cocycle} 
\phi(P\otimes U_g, Q\otimes U_{g^{-1}}) = \barint P U_g [\ln \Delta^{1/4}, Q U_{g^{-1}}] \end{equation*}
\end{prop}

\begin{pr} First, remark that 
\begin{align*}
\phi(P\otimes U_g, Q\otimes U_{g^{-1}}) &= \Pf_{z=0} \Tr \left([P U_g, Q U_{g^{-1}}]  \cdot \Delta^{-z/4}\right) \\
                          &= \Res_{z=0} \dfrac{1}{z} \Tr \left((P U_g Q U_{g^{-1}} - Q U_{g^{-1}}) P U_g  \cdot \Delta^{-z/4}\right)
\end{align*}

Then, work at $z \in \C$ with $\Re(z) \gg 0$, so that $\Tr \left([P U_g, Q U_{g^{-1}}]  \cdot \Delta^{-z/4}\right)$ is well-defined. Then, the trace property yields 
\begin{equation*}
\Tr \left([P U_g, Q U_{g^{-1}}]  \cdot \Delta^{-z/4}\right) = \Tr \left(P (U_g Q U_{g^{-1}}\Delta^{-z/4} -  U_g \Delta^{-z/4} Q U_{g^{-1}} )  \right)
\end{equation*}
Then, write 
\begin{align*}
U_g Q U_{g^{-1}}\Delta^{-z/4} - U_g \Delta^{-z/4} Q U_{g^{-1}}  &= U_g Q U_{g^{-1}}\Delta^{-z/4} - U_g Q \Delta^{-z/4} U_{g^{-1}} - U_g [\Delta^{-z/4},Q] U_{g^{-1}}  \\
                                 &= - U_g [\Delta^{-z/4},Q] U_{g^{-1}} + U_g Q U_{g^{-1}} [\Delta^{-z/4}, U_g] U_{g^{-1}} 
\end{align*}

To end the calculations, we need the following lemma, whose proof may be found in \cite{Hig2006}. 

\begin{lem} \label{CM trick} \emph{(Connes-Moscovici's trick, \cite{CM1995, Hig2006})} For every $z \in \C$, we have the following expansion,
\begin{align*}
& [\Delta^{-z}, Q] \sim \sum_{k\geq 1} \binom{-z}{k} Q^{(k)}\Delta^{-z-k} & U_{g^{-1}} [\Delta^{-z}, U_g] \sim \sum_{k\geq 1} \binom{-z}{k} U_{g^{-1}} U_g^{(k)} \Delta^{-z-k}
\end{align*}
where we denote $T^{(k)} = \ad(\Delta)^k(T)$, $\ad(\Delta) = [\Delta, \, . \,]$. 
\end{lem} 

Moreover, notice that for every integer $k \geq 1$, $Q^{(k)}$ and $U_{g^{-1}} U_g^{(k)}$ are classical Heisenberg pseudodifferential operators whose orders stricly decrease as $k$ grows. Hence, evaluating the expression under the trace, Theorem \ref{CMzeta} may be used. We deduce that on the one hand, the sums 
\begin{gather*}
\Res_{z=0} \dfrac{1}{z} \sum_{k\geq 1} \Tr \left(P \left(U_g \binom{-z/4}{k} Q^{(k)}\Delta^{-z/4 -k}  U_{g^{-1}} \right)\right) \\
\Res_{z=0} \dfrac{1}{z} \sum_{k\geq 1} \Tr \left(P \left(U_{g^{-1}} \binom{-z/4}{k} U_g^{(k)}\Delta^{-z/4-k}   \right)\right)
\end{gather*}
are finite, since the zeta function is holomorphic on a half-plane $\Re(z) \gg 0$. On the other hand, as the poles of the zeta function are simple, the terms carrying a power of $z^{2}$ vanish under the residue, and we are respectively left with 
\begin{gather*}
-\Res_{z=0} \sum_{k\geq 1} \Tr \left(P \left(U_g \dfrac{(-1)^{k-1}}{4k}   Q^{(k)}\Delta^{-z/4-k}  U_{g^{-1}} \right)\right) \\
-\Res_{z=0} \sum_{k\geq 1} \Tr \left(P \left(U_{g^{-1}} \dfrac{(-1)^{k-1}}{4k}  U_g^{(k)} \Delta^{-z/4-k}   \right)\right)
\end{gather*}
We then recognize the commutator with the logarithm of $\Delta^{1/4}$ (cf. \cite{Rod2013}), and we finally obtain 
\begin{align*}
\phi(P\otimes U_g, Q\otimes U_{g^{-1}}) &= \barint P \left(U_g [\ln \Delta^{1/4},Q] U_{g^{-1}} - U_g Q U_{g^{-1}} [\ln \Delta^{1/4}, U_g] U_{g^{-1}} \right) \\
                          &= \barint P U_g [\ln \Delta^{1/4}, Q U_{g^{-1}}]
\end{align*}
This ends the proof of the proposition. $\hfill{\square}$
\end{pr}

The pseudodifferential extension (\ref{extension}) is closely related to another extension. Indeed the quotient of $\Psi^0_{H,c}(M)$ by its two-sided ideal $\Psi_{H,c}^{-1}(M)$ of operators of order $\leq -1$ is $G$-equivariantly isomorphic to the commutative algebra of leading symbols $\cinfc(S^*_HM)$. The natural inclusion of smoothing operators in $\Psi_{H,c}^{-1}(M)$ and the leading symbol map thus yield a morphism of extensions
\begin{equation}
\vcenter{
\xymatrix{ 0 \ar[r] & \Psi_c^{-\infty}(M)\rtimes G \ar[r] \ar[d] & \Psi_{H,c}^0(M)\rtimes G \ar[r] \ar@{=}[d] & \S_{H,c}^0(M)\rtimes G \ar[r] \ar[d] & 0 \\
0 \ar[r] & \Psi_{H,c}^{-1}(M)\rtimes G \ar[r] & \Psi_{H,c}^0(M)\rtimes G \ar[r] & \cinfc(S_H^*M)\rtimes G \ar[r] & 0 }}
\end{equation}
The cyclic cohomology class of the operator trace $\Tr_{[1]}$ localized at unit extends in a straightforward manner to a cyclic cohomology class $[\tau]\in \HP^0(\Psi_{H,c}^{-1}(M)\rtimes G)$. The latter is represented, for any choice of even integer $k>v+2h$, by the cyclic $k$-cocycle $\tau_k$ defined as follows:
\begin{equation}
\tau_k(P_0,\ldots, P_k) = \Tr_{[1]} (P_0\ldots P_k)
\end{equation}
for all $P_i\in \Psi_{H,c}^{-1}(M)\rtimes G$. By naturality of excision, the class $\partial([\Tr_{[1]}]) \in \HP^1(\S_{H,c}^0(M)\rtimes G)$ is the pullback of $\partial ([\tau]) \in \HP^1(\cinfc(S_H^*M)\rtimes G)$ under the leading symbol homomorphism. Now the computation of $\partial([\tau])$ is fairly analogous to the above computation of $\partial([\Tr_{[1]}])$. We use the generalized Goodwillie theorem of Cuntz and Quillen \cite{CQ95}, which states that the periodic cyclic cohomology of an associative algebra $\A$ is isomorphic to the periodic cyclic cohomology of its \emph{completed tensor algebra}
\begin{equation}
\Th\A = \varprojlim_m T\A/(J\A)^m\ ,
\end{equation}
where the two-sided ideal $J\A \subset T\A$ is the kernel of the multiplication homomorphism $T\A\to\A$, $a_1\otimes\ldots\otimes a_k\mapsto a_1\ldots a_k$. We let $\A = \cinfc(S^*_HM)\rtimes G$ and choose any linear splitting 
\[\nu:\A\to \Psi^0_{H,c}(M)\rtimes G \] 
of the leading symbol homomorphism. $\nu$ is multiplicative up to the ideal $\Psi_{H,c}^{-1}(M)\rtimes G$. By the universal property of the tensor algebra, $\nu$ extends to an homomorphism $\nu_*:T\A \to \Psi^0_{H,c}(M)\rtimes G$ respecting the ideals, whence a morphism of extensions
\[
\xymatrix{ 0 \ar[r] & J\A \ar[r] \ar[d]_{\nu_*} & T\A \ar[r] \ar[d]_{\nu_*} & \A \ar[r] \ar@{=}[d] \ar@{.>}[dl]^{\nu} & 0 \\
0 \ar[r] & \Psi_{H,c}^{-1}(M)\rtimes G \ar[r] & \Psi_{H,c}^0(M)\rtimes G \ar[r] & \cinfc(S_H^*M)\rtimes G \ar[r] & 0 }
\]
Observe that the cocycle $\phi$, viewed as a cyclic 1-cocycle over $\Psi_{H,c}^0(M)\rtimes G$, vanishes on the large powers of the ideal $\Psi_{H,c}^{-1}(M)\rtimes G$ because it involves the Connes-Moscovici residue. Hence the composite $\phi\circ\nu_*$ extends to a cyclic 1-cocycle over $\Th\A$. By \cite{Per} Corollary 2.6, this cocycle is precisely a representative of the periodic cyclic cohomology class $\partial ([\tau]) \in \HP^1(\cinfc(S^*_HM)\rtimes G)$. Therefore we obtain

\begin{prop}\label{pcocycle}
Let $0\to \Psi_H^{-1}(M)\rtimes G \to \Psi_H^0(M)\rtimes G \to \cinf(S_H^*M)\rtimes G \to 0$ be the fundamental extension in the Heisenberg pseudodifferential calculus. Then the image $\partial([\tau])$ of the canonical trace localized at unit under the excision map $\partial : \HP^0(\Psi_{H,c}^{-1}(M)\rtimes G) \to \HP^1(\cinfc(S_H^*M)\rtimes G)$ is represented by a cyclic 1-cocycle over the completed tensor algebra of $\A=\cinfc(S_H^*M)\rtimes G$,
\begin{equation}
(\phi\circ\nu_*)(\ah_0,\ah_1) = \barint \big(\nu_*(\ah_0) [\ln \Delta^{1/4}, \nu_*(\ah_1)]\big)_{[1]} \ ,\qquad \forall \ \ah_0,\ah_1\in \Th\A\ ,
\end{equation}
for any choice of sub-Laplacian $\Delta$ and linear splitting $\nu:\A\to \Psi_{H,c}^0(M)\rtimes G$. 
\end{prop}

Then, we get a local index formula for any $G$-Heisenberg elliptic operator since as a Connes-Moscovici residue, it is an integral of the term of order $-v-2h$ in the asymptotic expansion of the operator $PU_g [\ln \Delta^{1/4}, QU_{g^{-1}}]$. However, extracting directly characteristic classes from it is in practice hard (except in low dimensions) because as a star-product, many higher derivatives of the symbols seem to be involved a priori. \\ 

The general strategy to handle the calculation is the construction of a cohomologous $(b,B)$-cocycle that depends only on certain symbols that already have degree $-v-2h$. This will reduce the cocycle at the level of the principal symbol and give an integral of a certain differential form over the Heisenberg cosphere bundle. This will be achieved through an "algebraic JLO formula" whose entries are formal Heisenberg symbols. This construction includes the following main steps: 
\begin{itemize}
\item[-] An algebra of operators acting on symbols,
\item[-] A related notion of "heat operator",
\item[-] A related notion of trace, 
\item[-] A related notion of "Dirac operator" $D$, in the sense that its commutator with a symbol gives (up to lower order terms) the differential of the symbol, 
\item Carrying all this to the equivariant setting. 
\end{itemize}
The constructions required are more tedious than the original JLO setting, however, they are much more flexible and has many advantages. As we will see, the constructions are purely algebraic, and we do not need to appeal to analytic properties of the heat equation, so that we can deal with much more general operators that Dirac type operators. In the same vein, we do not need to consider entire $(b,B)$-cochains. Thus, what we develop here holds in cases where Getzler's rescaling does not apply.

\section{Bimodule of Heisenberg formal symbols}\label{sbimod}

We adapt the framework developed in \cite{Per2012} to Heisenberg calculus on foliations. When the proof of a result is deferred to this paper, this means that it applies verbatim in our case. We shall mainly focus on the changes which occur in the Heisenberg setting. \\ 

Let $(M, V)$ be a foliated manifold of dimension $n$, where $V \subset TM$ is the sub-bundle of rank $v$ defining the foliation, and let $h$ denote the codimension of the foliation. We consider the $\Z_2$-graded algebra of formal Heisenberg symbols $\S_H(M,E)$, with $E=\Lambda^{\bullet}( T^* M\otimes\C)$. We view it as a left $\S_H(M,E)$-module and right $\PS_H(M,E)$-module : the left action of $a \in \S_H(M,E)$, and the right action $b \in \PS_H(M,E)$ on $\xi \in \S_H(M,E)$ are given by 
\[ a_L \cdot \xi = a\xi , \quad b_R \cdot \xi = \pm \xi b \]
Here, the sign $\pm$ depends on the parity of $b$ and $\xi$ : it is $-$ when both are odd and $+$ otherwise. This action defines a  $\Z_2$-graded subalgebra of $\End(\S_H(M,E))$ 
\[ \mathcal{L}(M) = \mathrm{span}\{a_L b_L \, ; \, a \in \S_H(M,E), b \in \PS_H(M,E) \} \]

Now, let us have a closer look on the operators contained in this algebra. Let $(x^1, \ldots , x^n)$ be a foliated coordinate system over an open subset $U\subset M$. The function $x^i$ is a symbol of order $0$, so that $x_L^i, x_R^i$ may be viewed as elements of $\mathcal{L}(M)$. The conjugate coordinate $p_i$ is a symbol of order $1$ if $i=1, \ldots, v$, but of order $2$ when $i=v+1, \ldots, n$. As before, this defines elements $p_{iL}, p_{iR}$ of $\mathcal{L}(M)$.  Now, observe that : 
\begin{align}
& (x_L^i - x_R^i) = \i \partial_{p_i}, \quad (p_{iL} - p_{iR}) = -\i \partial_{x^i} 
\end{align}
The same holds for the odd coordinates $\psi^{i}$ and $\overline{\psi}_i$ : 
\begin{align}
& (\psi_L^i - \psi_R^i) = \partial_{\overline{\psi}_i}, \quad (\overline{\psi}_{iL} - \overline{\psi}_{iR}) = \partial_{\psi^i} 
\end{align}
Hence, $\mathcal{L}(M)$ contains all the "elementary operations" on (Heisenberg) symbols. \\

Little calculations shows that a generic element $a_L b_R \in \L(M)$ reads over $U$ as a series
\begin{equation}
a_L b_R  = \sum_{\vert \alpha \vert = 1}^{k} \sum_{\vert \beta \vert = 1}^{\infty} \sum_{\vert \eta \vert = 1}^{n} \sum_{\vert \theta \vert = 1}^{n}  (s_{\alpha, \beta, \eta, \theta})_L (\psi^{\eta} \overline{\psi}^{\theta})_R \partial^{\alpha}_x \partial^{\beta}_p 
\end{equation} 
where $s_{\alpha, \beta, \eta, \theta} \in \S_H(U,E)$ and $k \in \N$. It is not necessarily true that any series of that form comes from an element of $\L(M)$.  \\

Now, consider the $\Z_2$-graded algebra $\mathscr{S}(M)=\L(M)[[\eps]]$ of \emph{formal power series} with coefficients $\L(M)$, and indeterminate $\eps$, which comes with a trivial grading. $\mathscr{S}$ is filtered by the subalgebras $\mathscr{S}_k(M) = \mathscr{S}(M)\eps^{k}$, for every $k\in \N$. This $k$ counts the minimal power of $\eps$ appearing in an element of $\mathscr{S}(M)$. We now define an important subalgebra of $\mathscr{S}(M)$. 

\begin{df} The subspace $\mathscr{D}^{m}(M) \subset \mathscr{S}(M)$ consists of elements $s = \sum s_k \eps^{k}$ such that in any {foliated} local chart $U \in M$, we have 
\begin{equation} \label{series D} s_k = \sum_{\vert \alpha \vert = 1}^{k} \sum_{\vert \beta \vert = 1}^{\infty} \sum_{\vert \eta \vert = 1}^{n} \sum_{\vert \theta \vert = 1}^{n}  (s^{m}_{k, \alpha, \beta, \eta, \theta})_L (\psi^{\eta} \overline{\psi}^{\theta})_R \partial^{\alpha}_x \partial^{\beta}_p  \end{equation}
where $s^{m}_{k, \alpha, \beta, \eta, \theta} \in \S_H(U,E)$ has \emph{Heisenberg order} $\leq m+(k+\vert \beta \vert - 3 \vert \alpha \vert) / 2$. We also denote $\mathscr{D}^{m}_k(M) = \mathscr{D}^{m}(M) \cap \mathscr{S}_k(M)$. We set 
\[ \mathscr{D}(M) = \bigcup_{m \in \R} \mathcal{D}^{m}(M) \]
\end{df}

The space $\mathscr{D}(M)$ is a bi-filtered subalgebra of $\mathscr{S}(M)$, that is $\mathscr{D}^{m}_k(M) \cdot \mathscr{D}^{m'}_{k'}(M) \subset \mathscr{D}^{m+m'}_{k+k'}(M)$, for $m, m' \in \R$ and $k, k' \in \N$ (\cite{Per2012}, Lemma 3.1). Using symbols \emph{with compact supports}, we define analogously the subalgebra $\mathscr{D}_c(M)\subset\mathscr{D}(M)$.\\ 

\begin{df} A \emph{generalized Laplacian} is an operator $\Delta \in \mathscr{D}^{1/2}_{1}$ of even parity, which can be written, in any local coordinate system over a local foliated chart $U \in M$ : 
\begin{equation} \label{generalized Laplacian}
\Delta \equiv \i \eps \partial_{x^{i}} \partial_{p_{i}} \text{mod } \mathscr{D}^{0}_{1}(U)
\end{equation} 
\end{df}
Throughout the paper, we shall use Einstein summation notation for repeated indices. \\

That such an operator exists is not obvious, cf. \cite{Per2012}, Lemma 3.3. We will see some important examples in the Section \ref{Dirac operators}. \\

A generalized Laplacian $\Delta$ will be our first point of departure towards the construction of a \emph{JLO formula on Heisenberg symbols}. In this type of formula, one needs to know how to deal with an exponential of such an operator in order to have a "heat operator". As a formal power series in $\eps$ this indeed defines an element of $\mathscr{S}(M)$ : 
\begin{equation}
\exp(t\Delta) = \sum_{k=0}^{\infty} \dfrac{t^{k}}{k!} \Delta^{k}, \, \forall t \in \R
\end{equation}
This operator does not belong to $\mathscr{D}(M)$. We define a one parameter group of automorphisms $(\sigma_\Delta^{t})_{t \in \R}$ of the algebra $\mathscr{S}(M)$ as follows : 
\begin{equation}
\sigma_\Delta^{t}(s) = \exp(t\Delta) s \exp(-t\Delta), \, \forall s \in \mathscr{D}(M)
\end{equation}

\begin{lem} For every generalized Laplacian $\Delta$, $(\sigma_\Delta^{t})_{t \in \R}$ is actually a one parameter group of automorphisms of $\mathscr{D}(M)$. More precisely, one has, for every $m \in \R$, $k \in \N$ : 
\begin{equation*}
[\Delta, \mathscr{D}^{m}_k] \subset \mathscr{D}^{m}_{k+1}, \quad \sigma_\Delta^{t}(\mathscr{D}^{m}_k) \subset \mathscr{D}^{m}_k
\end{equation*}
\end{lem}

\begin{prop} \textup{\emph{(Duhamel formula)}} \label{(Duhamel formula)} Let $\Delta + s$ be a perturbation of a generalized Laplacian $\Delta$, where $s \in \mathscr{D}_1^{0}(M)$. Then, 
\begin{equation}
\exp(\Delta + s) = \sum_{k=0}^{\infty} \int_{\Delta_k} \exp(t_0 \Delta) s \exp(t_1 \Delta) \ldots s \exp(t_k \Delta) \, dt
\end{equation}
where $\Delta_k$ is the standard $k$-simplex, and $dt = dt_1 \ldots dt_k$. Equivalently,
\begin{equation}
\exp(\Delta + s) = \sum_{k=0}^{\infty} \int_{\Delta_k} \sigma^{t_0}_\Delta(s) \sigma^{t_0 + t_1}_\Delta(s) \ldots \sigma^{t_0 + \ldots t_{k-1}}_\Delta(s) \exp(\Delta) \, dt
\end{equation}
\end{prop}

\begin{df} Let $\Delta$ be a generalized Laplacian.  The \emph{bimodule of trace class operators} is the subspace $\mathscr{T}(M) = \mathscr{D}_c(M)\exp(\Delta)$ of $\mathscr{S}(M)$.
\end{df}

\begin{rk} $\mathscr{T}(M)$ is a $\mathscr{D}(M)$-bimodule, and does not depend on the choice of the generalized Laplacian (\cite{Per2012}, Proposition 3.6). Notice that $\mathscr{T}(M)$ is not closed under multiplication. 
\end{rk} 

The terminology will be explained in the next section. Finally notice that if $G\subset\Diff(M)$ is a group of foliated diffeomorphisms, the algebra $\D(M)$ and the bimodule $\mathscr{T}(M)$ carry natural $G$-actions by automorphisms.

\section{Canonical trace on the bimodule $\mathscr{T}$}\label{scano}

\subsection{Construction of the trace}
Let $(M^{n}, V)$ be a foliated manifold of codimension $h$, and let $v$ denote the dimension of the leaves, so that $n=v+h$, and take $E=\Lambda^{\bullet} (T^{*}M\otimes\C)$. The aim of this section is to construct a canonical trace on $\mathscr{T}(M)$ from the Wodzicki residue \ref{CM residue}. \\ 

First, we work locally. Let $U \subset M$ be a foliated local chart of $M$ . Recall that $\mathscr{T}(U)$ is a bimodule over $\mathscr{D}(U)$. A trace on $\mathscr{T}(U)$ is in this sense a linear map $\mathscr{T}(U) \to \C$ vanishing on the subspace $[\mathscr{T}(U), \mathscr{D}(U)]$ of graded commutators. Choose a coordinate system $(x,p)$ on $T^{*}U$ adapted to the foliation, and let $\Delta$ be the "flat" (generalized) Laplacian  
\[ \Delta = \i \eps \partial_{x^{i}} \partial_{p_{i}} \]

For every multi-indices $\alpha$ and $\beta$, we define a bracket operation 
\begin{equation}
\langle \partial_x^{\alpha} \partial_p^{\beta} \exp \Delta \rangle = \partial_x^{\alpha} \partial_p^{\beta} \exp \left. \left(\dfrac{\i}{\eps} (p_i - q_i) (x^{i} - y^{i}) \right) \right\vert_{x=y, p=q}
\end{equation}
Remark that this vanishes unless $\vert \alpha \vert = \vert \beta \vert$. \\

\begin{ex} One has 
\[ \langle \exp \Delta \rangle =1,\quad \langle \partial_{x^{i}}\exp \Delta \rangle = \langle \partial_{p_{i}}\exp \Delta \rangle = 0, \quad \langle \partial_{x^{i}}\partial_{p_{j}}\exp \Delta \rangle = \dfrac{\i}{\eps} \delta_i^{j}  \] 
where $\delta_i^{j}$ denotes the Kronecker symbol. More generally, the formula with a polynomial $\partial_{x}^{\alpha}\partial_{p}^{\beta}$ involves all possible contractions between $\partial_{x^{i}}$ and $\partial_{p_{j}}$. For example,
\[ \langle \partial_{x^{i}}\partial_{x^{j}}\partial_{p_{k}}\partial_{p_{l}} \exp \Delta \rangle = \left(\dfrac{\i}{\eps}\right)^{2}(\delta_i^{k}\delta_j^{l} + \delta_i^{l}\delta_j^{k}) \]
\end{ex}

We also define a contraction map for the odd variables : for every multi-indices $\eta$ and $\theta$, we set
\begin{equation}
\langle (\psi^{\eta} \overline{\psi}^{\theta})_R \rangle = (-1)^{n} \tr_s(\psi^{\eta} \overline{\psi}^{\theta})
\end{equation}
where $\tr_s$ is the graded trace of endomorphisms of the vector bundle $E$. In particular we have $\langle (\psi^{1} \ldots \psi^{n} \overline{\psi}_{n} \ldots \overline{\psi}_{1})_R \rangle = 1$

From this, we construct a linear map 
\begin{equation}
\langle\langle \ .\ \rangle\rangle : \mathscr{T}(U) \to \S_H(U,E)[[\eps]]
\end{equation}
as follows. Let $s \exp(\Delta) \in \mathscr{T}(U)$ be a generic element, where $s = \sum_{k \geq 0} s_k \eps^{k} \in \mathscr{D}^{m}(U)$. The symbol $s_k$ may be written 
\[s_k = \sum_{\vert \alpha \vert = 1}^{k} \sum_{\vert \beta \vert = 1}^{\infty} \sum_{\vert \eta \vert = 1}^{n} \sum_{\vert \theta \vert = 1}^{n}  (s_{k, \alpha, \beta, \eta, \theta})_L (\psi^{\eta} \overline{\psi}^{\theta})_R \partial^{\alpha}_x \partial^{\beta}_p \]
where $s_{k, \alpha, \beta, \eta, \theta} \in \S_H(U,E)$ has \emph{Heisenberg order} $\leq m+(k+\vert \beta \vert - 3 \vert \alpha \vert) / 2$. Then, we set

\[ \langle\langle s_k \exp \Delta \rangle\rangle = \sum_{\vert \alpha \vert = 1}^{k} \sum_{\vert \beta \vert = 1}^{\infty} \sum_{\vert \eta \vert = 1}^{n} \sum_{\vert \theta \vert = 1}^{n}  s_{k, \alpha, \beta, \eta, \theta} \langle(\psi^{\eta} \overline{\psi}^{\theta})_R \rangle \langle \partial^{\alpha}_x \partial^{\beta}_p \exp \Delta \rangle \]

This sum is finite by definition of the even contraction, and is consequently a polynomial of degree at most $k$ in the variable $\eps^{-1}$. Then, $\langle\langle s_k \exp \Delta \rangle\rangle \eps^{k}$ is a polynomial of degree at most $k$ in $\eps$. Then, we finally define 
\begin{equation}
\langle\langle s \exp \Delta \rangle\rangle = \sum_{k \geq 0} \langle\langle s_k \exp \Delta \rangle\rangle \eps^{k}
\end{equation}  
which is an element of $\S_H(U,E)[[\eps]]$ (this is not totally obvious, refer to \cite{Per2012}, Lemma 4.1). 

We can now pass to the definition of the trace on $\mathscr{T}(M)$. 

\begin{df} Let $U \in M$ be a foliated local chart, and denote by $\langle\langle s \exp \Delta \rangle\rangle[n]\in S_H(U,E)$ the coefficient of $\eps^{n}$ in the formal series $\langle\langle s \exp \Delta \rangle\rangle$. Then, we define the following graded trace
\begin{equation}
\Tr^{U}_s : \mathscr{T}(U) \to \C, \quad  \Tr^{U}_s(s \exp \Delta) = \barint \langle\langle s \exp \Delta \rangle\rangle[n]
\end{equation} 
This map does not depend on the choice of foliated coordinates $(x,p)$ on $T^{*}U$, so that these maps may be glued together to give a canonical graded trace :
\begin{equation}
\Tr_s : \mathscr{T}(M) \to \C
\end{equation}
on the $\mathscr{D}(M)$-bimodule of trace class operators. 
\end{df}

The proof that this is a trace is the same as that of \cite{Per2012}, Lemma 4.2. That we can glue these quantities to get a global functional on the whole foliation $M$ is Proposition 4.3 of the same paper. For the same reason, $\Tr_s$ is invariant under the action of any group $G$ of foliated diffeomorphisms on $\mathscr{T}(M)$. 
 
\subsection{An algebraic Mehler formula}

In this section, we recall how the "Todd series" can be recovered from the contractions we defined in the previous paragraph (\cite{Per2012}). The formula may be seen as a pseudodifferential analogue of the Mehler formula for the harmonic oscillator, and will be crucial for obtaining the index theorem. We keep the notations of the previous subsection and work in the foliated local chart $U$. \\

For a $N\times N$ matrix $\Omega$ with coefficients in $\C[[\eps]]$, which has no degree zero term in $\eps$, we can define the following formal power series in $M_N(\C[[\eps]])$ and in $\C[[\eps]]$ :
\begin{equation}
\dfrac{\Omega}{e^{\Omega} - 1} = 1 - \dfrac{1}{2}\Omega + \dfrac{1}{12}\Omega^{2} + \ldots \, , \quad \Td(\Omega) = \det \left(\dfrac{\Omega}{e^{\Omega} - 1} \right)
\end{equation}  
$\Td(\Omega)$ is the \emph{Todd series} of $\Omega$. \\

Now, consider the operator 
\[ p_L \cdot \Omega \cdot \partial_p =  p_{iL} \cdot \Omega^{i}_j \cdot \partial_{p_{j}} \]  
and the perturbation of the flat Laplacian $\Delta + p_L \cdot \Omega \cdot \partial_p$, which is no longer a generalized Laplacian. However, by the Duhamel formula, Proposition \ref{(Duhamel formula)}, $\exp (\Delta + p_L \cdot \Omega \cdot \partial_p)$ still defines an element of $\mathscr{T}(U)$.  

\begin{prop} For every multi-indices $\alpha$ and $\beta$, we have 
\begin{equation}
\langle \partial_x^{\alpha} \partial_p^{\beta} \exp (\Delta + p_L \cdot \Omega \cdot \partial_p)\rangle = \Td(\Omega)s(\Omega,p) 
\end{equation}
where the symbol $s(\Omega,p)$ is polynomial in $p$ and given by the following formula :
\[ s(\Omega,p) = \left. \partial_x^{\alpha} \partial_p^{\beta} \exp \left(\dfrac{\i}{\eps}q \cdot \Omega \cdot (x-y) + \dfrac{\i}{\eps}(p-q) \cdot \dfrac{\Omega}{1-e^{-\Omega}} \cdot (x-y) \right)\right\vert_{x=y,p=q}  \]
\end{prop}

\begin{ex} \label{todd matrix} We give some particular cases of the formula which will be useful in the sequel. We have
\begin{gather}
\langle \exp (\Delta + p_L \cdot \Omega \cdot \partial_p) \rangle = \Td(\Omega) \\
\langle (\i \eps \partial_x + p_L \cdot \Omega)^{\alpha}  \exp(\Delta + p_L \cdot \Omega \cdot \partial_p) \rangle = 0
\end{gather}
where $\alpha$ is any non-zero multi-index. 
\end{ex} 

\section{Dirac operators}\label{sdirac}

\subsection{Generalities}

The algebra of differential forms $\Omega^{\bullet}(M)$ on $M$ and the Lie algebra $\Vect(M)$ of vector fields  may be seen as elements of the space $\mathscr{PS}_H^{0}(M,E)$ of polynomials Heisenberg symbols of degree $0$ via the following maps : 
\begin{align*}
&\mu : dx^{i_1} \wedge \ldots \wedge dx^{i_k} \in \Omega^{\bullet}(M) \mapsto \psi^{i_1} \ldots \psi^{i_k} \in \mathscr{PS}_H^{0}(M,E), \\
&\iota : \partial_{x^{i}} \in \Vect(M) \mapsto \iota(\partial_{x^{i}}) = \overline{\psi}_i \in \mathscr{PS}_H^{0}(M,E)
\end{align*}

Then, we shall be interested in various subspaces of $\L(M) = \S_H(M,E)_L \mathscr{PS}_H(M,E)_R$, which are needed to see where lie the generalized Dirac operators. \\ 

Let $\mathcal{SPS}_H^{1}(M,E) \subset \mathscr{PS}_H^{1}(M,E)$ be the space of differential operators $a$ of order 1, with \emph{scalar Heisenberg principal symbol}, whose local expression reads 
\begin{equation*}
 a(x,p) = a^{i}(x)p_i + a^{i}_j(x) \psi^{j} \overline{\psi}_i + b(x) 
\end{equation*} 
where the coefficients depending on $x$ are smooth functions. Remark that $\mathcal{SPS}_H^{1}(M,E)$ is a Lie algebra. We also consider the subspaces $\mathcal{SPS}_H^{1}(M,E)_L \Omega^{1}(M)_R \subset \mathscr{D}_0^{1}(M)$ and $\Omega^{0}(M)_L \Vect(M)_R \subset \mathscr{D}_0^{0}(M)$, respectively spanned by elements which are locally given by the series
\begin{gather*}
s = \sum_{\vert \alpha \vert \geq 0} (s_{\alpha i}^{k}(x) p_{k} + s_{\alpha i j}^{k}(x) \psi^{i} \overline{\psi}_{k} + s_{\alpha i}(x))_L \psi^{i}_R \partial_p^{\alpha} \\
 r = \sum_{\vert \alpha \vert \geq 0} (r_{\alpha}^{i})_L \overline{\psi}_{iR} \partial_{p}^{\alpha} 
\end{gather*}
The coefficients depending on $x$ are again smooth functions. 

\begin{df} \label{Dirac operators} A \emph{generalized Dirac operator} $D$ is an element of $\mathscr{D}(M)$ of the form 
\begin{equation*}
D = i \eps \nabla + \overline{\nabla} \in \mathscr{D}_1^{1}(M) + \mathscr{D}_0^{-1/2}(M)
\end{equation*}
where $\nabla$ and $\overline{\nabla}$ are such that
\begin{align*}
&\nabla \equiv \psi^{i}_{R} \partial_{x^{i}} \, \text{mod } \mathcal{SPS}_H^{1}(M,E)_L \Omega^{1}(M)_R, \\
&\overline{\nabla} \equiv \overline{\psi}_{iR} \partial_{p_{i}} \text{mod } \Omega^{0}(M)_L \Vect(M)_R \cap \mathscr{D}_0^{-1}(M)
\end{align*}
\end{df}

Remark that if $G$ is a group of foliated diffeomorphisms on $M$, it transforms Dirac operators into Dirac operators. The terminology lies in the following important proposition (\cite{Per2012}). 

\begin{prop} Let $D$ be a generalized Dirac operator. Then, $-D^{2}$ is a generalized Laplacian. 
\end{prop}  

The rest of the paragraph gives the two crucial examples of generalized Dirac operators. The proof of the formulas can be found in \cite{Per2012}.

\subsection{De Rham - Dirac operators} The exterior differentiation $d$ acting on the space of differential forms $\Omega^{\bullet}(M)$ defines an element of $\mathcal{PS}^{1}(M,E)$. Its right action on the bimodule of formal Heisenberg symbols $\S_H(M,E)$ gives an element of odd degree $d_R \in \L(M)$. Locally, 
\begin{equation}
 d_R = \i (p_i \psi^{i})_R = \i \psi^{i}_R p_{iR} = -\psi^{i}_R \partial_{x^{i}} + \i p_{iL} \psi^{i}_R
\end{equation} 
As $\i p_{iL} \psi^{i}_R \in \mathcal{SPS}^{1}(M,E)$, we have a generalized Dirac operator 
\begin{equation}
D = -\i \eps d_R + \overline{\nabla}
\end{equation}
for any choice of $\overline{\nabla}$ as in Definition \ref{Dirac operators}. Such generalized Dirac operators will be called of \emph{de Rham - Dirac type}. 

\begin{prop} Let $D = -\i \eps d_R + \overline{\nabla}$ be a de Rham - Dirac operator. Locally, the associated generalized Laplacian is given by the following formula : 
\begin{multline}
-D^{2} = \i \eps \left( \partial_{x^{i}} \partial_{p_{i}} + \sum_{\vert \alpha \vert \geq 2} (a_{\alpha}^{i})_L \partial_{x^{i}} \partial^{\alpha}_{p} \right) + \eps \left( p_{iL} \partial_{p_{i}} + \sum_{\vert \alpha \vert \geq 2} (a_{\alpha}^{i} p_i)_L  \partial^{\alpha}_{p} \right) \\
+ \eps \left( (\psi^{i} \overline{\psi}_i)_R + \sum_{\vert \alpha \vert \geq 1} (b_{\alpha j}^{i})_L (\psi^{j} \overline{\psi}_i)_R \partial^{\alpha}_{p} \right)
\end{multline}
where the coefficients $a_{\alpha}^{i}, b_{\alpha j}^{i}$ are smooth functions.
\end{prop}

\subsection{Dirac operators associated to affine connections}

Let $\Gamma$ be a torsion-free affine connection on $M$, characterized by its Christoffel symbols $\Gamma_{ij}^{k}$ in a local coordinate system $(x^1,\ldots,x^n)$ over $U$. Then, we define a "covariant derivative" operator on $\S_H(U,E)$ given by 
\begin{equation}
\nabla_i^{\Gamma} = \partial_{x^{i}} + \Gamma_{ij}^{k}(x)_L\left( p_{kL} \partial_{p_j} + (\overline{\psi}_k \psi^{j})_L - (\overline{\psi}_k \psi^{j})_R \right) 
\end{equation} 
This is not properly speaking a covariant derivative, since the coordinates $x$ and $p$ do not commute. However, the action of $\nabla_i^{\Gamma}$ on the generators $x,p,\psi,\overline{\psi}$ are what we expect from a covariant derivative : 
\begin{align}
&\nabla_i^{\Gamma}(x^{k}) = \delta^{k}_{i}, && \nabla_i^{\Gamma}(p_j) =  \Gamma_{ij}^{k} p_k, && \nabla_i^{\Gamma}(\psi^{k}) = - \Gamma_{ij}^{k} \psi^{j}, && \nabla_i^{\Gamma}(\overline{\psi}_{j}) = \Gamma_{ij}^{k} \overline{\psi}_{k}
\end{align}

A generalized Dirac operator $D = \i \eps \nabla + \overline{\nabla}$ is called \emph{affiliated to the affine connection} $\Gamma$ on $M$ if locally over $U$ we have
\begin{equation}
\nabla = \psi^{i}_R(\nabla_i^{\Gamma} + s)
\end{equation}
with $s \in \mathcal{SPS}_H^{1}(M,E)_L \Omega^{1}(M)_R \cap \mathscr{D}_0^{0}(M)$. 

\begin{prop} \label{lichner} For such a Dirac operator $D$, one has the following analogue of the Lichnerowicz formula : 
\begin{multline}
-D^{2} = \i \eps \left( \partial_{x^{i}} \partial_{p_{i}} + (\Gamma_{ij}^{k})_L (\psi^{i}\overline{\psi}_k)_R \partial_{p_{j}} + u + v \right) \\ 
+ \eps^{2}\left(\dfrac{1}{2} (\psi^{i} \psi^{j})_R (R_{lij}^{k})_L (p_{kL} \partial_{p_{l}} + (\overline{\psi}_k \psi^{l})_L) + w \right)
\end{multline}
where $R_{lij}^{k} = \partial_{x^{i}} \Gamma_{jl}^{k} - \partial_{x^{j}} \Gamma_{il}^{k} + \Gamma_{im}^{k}\Gamma_{jl}^{m} - \Gamma_{jm}^{k}\Gamma_{il}^{m}$ are the components of the curvature tensor of $\Gamma$, and
\begin{gather*}
u = \sum_{\vert \alpha \vert \geq 2} \left( (u_{\alpha i})_L \partial_{x^{i}} + (u_{\alpha}^{k} p_k)_L + (u_{\alpha i}^{k})_L (\psi^{i}\overline{\psi}_k)_R + (u_{\alpha})_L \right) \partial^{\alpha}_p \\
v = \sum_{\vert \alpha \vert \geq 1} (v_{\alpha i}^{k} \overline{\psi}_k \psi^i)_L \partial^{\alpha}_p \\
w = (\psi^{i} \psi^{j})_R \left(\sum_{\vert \alpha \vert \geq 2}(w_{\alpha i j}^{k} p_k)_L \partial^{\alpha}_p + \sum_{\vert \alpha \vert \geq 1}(w_{\alpha l i j}^{k} \overline{\psi}_k \psi^{l} + w_{\alpha i j})_L \partial^{\alpha}_p \right)
\end{gather*}
where  the coefficients are smooth functions on $M$. 
\end{prop}  

\section{Equivariant cohomology} \label{sequiv}

Let $G$ be a discrete group acting by orientation-preserving diffeomorphisms on a smooth oriented manifold $M$. Following \cite{Per2014}, we shall explain an alternative  construction of Connes' characteristic map from the $G$-equivariant cohomology of $M$ to the periodic cyclic cohomology of the crossed product algebra $\cinfc(M)\rtimes G$. Our approach differs slightly from the original construction given by Connes in \cite{C83} and looks quite complicated at first sight, but it is particularly well-adapted to the proof of the equivariant index theorem. \\

\textbf{Note. In this section only, the index $n$ does not denote the dimension of the manifold $M$. 
}
\subsection{Classifying spaces}

We recall that the nerve of a discrete group $G$ is the simplicial set $NG_\bullet$ with $NG_n = G^n$ for all $n\geq 0$. The face maps $\delta_i:NG_n\to NG_{n-1}$ and degeneracy maps $\sigma_i: NG_n\to NG_{n+1}$ are given by
\begin{multline}
\delta_0(g_1,\ldots,g_n) = (g_2,\ldots,g_n)  \\
\delta_i(g_1,\ldots,g_n) = (g_1,\ldots,g_ig_{i+1},\ldots,g_n)\qquad 1\leq i\leq n-1  \\
\delta_n(g_1,\ldots,g_n) = (g_1,\ldots,g_{n-1})  \\
\sigma_i(g_1,\ldots,g_n) = (g_1,\ldots,g_i,1,g_{i+1},\ldots,g_n) \qquad 0\leq i\leq n \ .\\
\end{multline}
Let $\Delta_n=\{(s_0,\ldots,s_n)\in [0,1]^{n+1}\ |\ s_0+\ldots+s_n=1\}$ be the standard $n$-simplex in $\R^{n+1}$, with $\delta^i:\Delta_{n}\to\Delta_{n+1}$, $(s_0,\ldots,s_n)\mapsto (s_0,\ldots,s_{i-1},0,s_{i},\ldots,s_n)$ the inclusion of the $i$-th face, and $\sigma^i:\Delta_{n}\to \Delta_{n-1}$, $(s_0,\ldots,s_n)\mapsto (s_0,\ldots,s_i+s_{i+1},\ldots,s_n)$ the collapse of the $i$-th edge. The classifying space of $G$ is the geometric realization of the simplicial set $NG_\bullet$, defined as the quotient
\begin{equation}
BG = \Big(\bigcup_{n\geq 0} NG_n\times \Delta_n\Big)/\sim
\end{equation}
where the equivalence relation $\sim$ identifies a point $(g,\delta^i s)\in NG_n\times \Delta_n$ (resp. $(g,\sigma^i s)\in NG_n\times \Delta_n$) with the point $(\delta_ig,s)\in NG_{n-1}\times \Delta_{n-1}$ (resp. $(\sigma_ig,s)\in NG_{n+1}\times \Delta_{n+1}$). Let $\Omega(\Delta_n)$ denote the Differential Graded (DG) algebra of complex-valued smooth differential forms over $\Delta_n$ which are extendable over the hyperplane $\{(s_0,\ldots,s_n)\in \R^{n+1}\ |\ s_0+\ldots+s_n=1\}$. Let $\Omega(NG_m\times \Delta_n)$ be the space of functions from the discrete set $NG_m$ to $\Omega(\Delta_n)$. This is naturally a DG algebra. A differential form $\omega$ of degree $k$ over $BG$ is a collection of $k$-forms $\omega_n\in \Omega^k(NG_n\times\Delta_n)$, $n\in\N$, subject to the constraints
\[
(\Id\times\delta^i)^*\omega_n = (\delta_i\times\Id)^*\omega_{n-1}\ ,\qquad (\Id\times\sigma^i)^*\omega_n = (\sigma_i\times\Id)^*\omega_{n+1}\ ,
\]
for all $i=0,\ldots, n$ and $n\geq 0$. The space $\Omega(BG)$ of differential forms over $BG$ is a DG algebra. The de Rham cohomology of $BG$, which is the cohomology of the complex $\Omega(BG)$ endowed with the exterior differential, is known to be canonically isomorphic to the group cohomology of $G$ with complex coefficients:
\begin{equation}
H^\bullet(BG) \cong H^\bullet(G,\C)\ .
\end{equation}
The universal $G$-bundle over the nerve $NG_\bullet$ is the simplicial set $N\overline{G}_\bullet$ with $N\overline{G}_n = G^{n+1}$ for all $n$, and the face and degeneracy maps are
\begin{multline}
\delta_i(g_0,\ldots,g_n) = (g_0,\ldots,\check{g}_i,\ldots,g_n)\qquad 0\leq i\leq n  \\
\sigma_i(g_0,\ldots,g_n) = (g_0,\ldots,g_i,g_i,\ldots,g_n) \qquad 0\leq i\leq n \ ,\\
\end{multline}
where the symbol $\check{}$ denotes omission. The projection $N\overline{G}_\bullet \to NG_\bullet$ defined by $(g_0,\ldots,g_n)\mapsto (g_0g_1^{-1},g_1g_2^{-1},\ldots ,g_{n-1}g_n^{-1})$ is a simplicial map. Its fibers are in one to one correspondence with the orbits of the (free) $G$-action $(g_0,\ldots,g_n)\cdot g = (g_0g,\ldots,g_ng)$, which is also a simplicial map for all $g\in G$. The geometric realization 
\begin{equation}
EG = \Big(\bigcup_{n\geq 0} N\overline{G}_n\times \Delta_n\Big) /\sim
\end{equation}
is therefore a $G$-bundle over $BG$. The DG algebra of differential forms $\Omega(EG)$ defined as above carries a natural action of $G$. By pullback, $\Omega(BG)$ is isomorphic to the DG subalgebra of $\Omega(EG)$ consisting of $G$-invariant differential forms. Hence $H^\bullet(BG)$ is also the cohomology of the complex of $G$-invariant differential forms on $EG$.\\

Now let $G\subset\Diff(M)$ be a discrete group acting by diffeomorphisms on a smooth manifold $M$. The product $EG\times M$, endowed with the diagonal $G$-action, is a $G$-bundle over the quotient $EG\times_G M$. The algebra $\Omega(EG\times M)$, defined as the collection of differential forms over the manifold $N\overline{G}_n\times\Delta_n\times M$ with gluing constraints as above, is naturally endowed with the total differential coming from $EG$ and $M$, and inherits an action of $G$ by pullback. The DG subalgebra of $G$-invariant differential forms is isomorphic to $\Omega(EG\times_G M)$. We define the $G$-equivariant cohomology of $M$ (with complex coefficients) as the corresponding de Rham cohomology $H^\bullet(EG\times_G M)$.\\  

The Chern-Weil theory of characteristic classes for vector bundles is easily generalized to the equivariant case. Let $V$ be a $G$-equivariant (complex) vector bundle over $M$, and choose a connection $\Gamma_0$ on $V$. Of course $\Gamma_0$ cannot be chosen invariant for general actions of $G$, and we denote by $\Ad_g(\Gamma_0)$ its image under conjugation by an element $g\in G$. The set of all connections being an affine space, at any point $(g_0,\ldots,g_n)(s_0,\ldots,s_n)\in N\overline{G}_n\times\Delta_n$ we can build a new connection $\Gamma$ on $V$ by means of the barycentric formula
\begin{equation}
\Gamma(g_0,\ldots,g_n)(s_0,\ldots,s_n) = \sum_{i=0}^n s_i\, \Ad_{g_i}^{-1}(\Gamma_0)\ .
\end{equation}
Let $W=EG\times V$ be the pullback of the vector bundle $V$ over $EG\times M$. If $d_H$ denotes the exterior differential over $EG$, then $d_H+\Gamma$ is a connection on $W$, whose curvature 2-form $R = [d_H,\Gamma]+\Gamma^2 \in \Omega^2(EG\times M, \End(W))$ reads
\[
R(g_0,\ldots,g_n)(s_0,\ldots,s_n) = \sum_{i=0}^n ds_i\, \Ad_{g_i}^{-1}(\Gamma_0) + \sum_{i,j} s_is_j\, \Ad_{g_i}^{-1}(\Gamma_0)\, \Ad_{g_j}^{-1}(\Gamma_0)\ .
\]
Observe that $R$ is always the sum of a form of bidegree $(1,1)$ and a form of bidegree $(0,2)$ with respect to the product manifold $EG\times M$. Since $R$ is $G$-equivariant by construction, any $\Ad$-invariant polynomial in the curvature yields a closed $G$-invariant differential form on $EG\times M$. In particular the Chern character $\ch(V)$ and the Todd class $\Td(V)$ are represented by 
\begin{equation}
\ch(R) = \tr(\exp(R))\ ,\qquad \Td(R) = \det\left(\frac{R}{e^{R}-1} \right)\ ,
\end{equation}
and a classical homotopy argument shows that their respective cohomology classes in $H^\bullet(EG\times_G M)$ do not depend on the particular choice of connection $\Gamma_0$ on $V$. Notice that we do not insert the standard $\i/2\pi$ factors in front of the curvature in these formulas, because the normalization we choose here is more natural with respect to the computation of index pairings using cyclic cohomology. \\

\subsection{Characteristic map}

We now explain our construction of Connes' characteristic map from the $G$-equivariant cohomology of $M$ to the periodic cyclic cohomology of $\A=\cinfc(M)\rtimes G$. The idea is to twist the universal tensor extension of the group ring $\C G$
\[
0 \to J\C G \to T\C G \to \C G \to 0
\]
by the DG algebra of smooth differential forms on $EG\times M$. Indeed $G$ acts on both spaces $EG$ and $M$ (from the right), and the induced action (from the left) by pullback on the graded-commutative algebra of differential forms $\Omega(EG\times M)$ commutes with the total differential $d$. Let $\Omega_p(EG\times M)$ be the subalgebra of differential forms $\alpha \in \Omega(EG\times M)$ which have compact $M$-support at any point of $EG$ (the subscript $_p$ stands for "proper"). The crossed product
\begin{equation}
\G = \Omega_p(EG\times M)\rtimes G
\end{equation}
is naturally a DG algebra. The product of two elements reads
\[ 
(\alpha\otimes U_{g_1})(\beta\otimes U_{g_2}) = \alpha \wedge U_{g_1}(\beta) \otimes U_{g_1g_2}
\]
for all $\alpha,\beta \in \Omega_p(EG\times M)$ and $g_i \in G$, where $U_{g_1}(\beta)$ is the pullback of $\beta$ by the diffeomorphism $g_1$. The differential reads $d(\alpha\otimes U_g)=d\alpha\otimes U_g$. An algebra extension of $\G$ is defined as the vector space
\begin{equation}
\H = \Omega_p(EG\times M)\otimes T\C G\ ,
\end{equation}
graded by the differential form degree ($T\C G$ is trivially graded), and endowed with the twisted product 
\[
(\alpha\otimes U_{g_1}\otimes\ldots\otimes U_{g_n})(\beta\otimes U_{g_{n+1}}\otimes\ldots\otimes U_{g_{n+m}}) = \alpha\wedge U_{g_1\ldots g_n}(\beta) \otimes U_{g_1}\otimes\ldots\otimes U_{g_{n+m}}\ .
\]
Also the de Rham differential is extended to $\H$ by $d(\alpha\otimes U_{g_1}\otimes\ldots\otimes U_{g_n}) = d\alpha\otimes U_{g_1}\otimes\ldots\otimes U_{g_n}$. Clearly the obvious multiplication map $\H\to\G$, $\alpha\otimes U_{g_1}\otimes\ldots\otimes U_{g_n} \mapsto \alpha\otimes U_{g_1 \ldots g_n}$ is a morphism of DG algebras, hence its kernel $\I = \Omega_p(EG\times M)\otimes J\C G$ is a two-sided DG ideal. We define a completion of $\H$ as  
\begin{equation}
\Hh = \bigoplus_{k\geq 0}\varprojlim_n \big(\Omega_p^k(EG\times M)\otimes T\C G/(J\C G)^n\big) \ .\label{toto}
\end{equation}
The product and differential on $\H$ extend in an obvious way to $\Hh$. If $EG$ were a finite-dimensional manifold, the sum over the form degree $k$ would be finite and $\Hh=\varprojlim_n \Omega_p(EG\times M)\otimes T\C G/(J\C G)^n$ would coincide with the $\I$-adic completion of $\H$ as in \cite{Per2014}. This however does not hold for the construction of $EG$ above, and (\ref{toto}) is really a smaller algebra.\\

We view $\cinfc(M)\subset\Omega^0_p(EG\times M)$ as the subalgebra of 0-forms which are constant in the direction $EG$. This identification is $G$-equivariant and extends to a morphism of algebras $\rho:\A\to\G$. The universal property of the tensor algebra thus yields an homomorphism $\rho_*: T\A \to \H$ explicitly given by 
\[
\rho_*(f_1 U_{g_1}\otimes f_2 U_{g_2}\ldots\otimes f_n U_{g_n}) = f_1 U_{g_1}(f_2)\ldots U_{g_1\ldots g_{n-1}}(f_n) \otimes U_{g_1}\otimes U_{g_2} \otimes\ldots\otimes U_{g_n}
\]
on any $n$-tensor, $f_i\in\cinfc(M)$, $g_i\in G$. One easily checks that $\rho_*$ carries the ideal $(J\A)^n$ to $\Omega_p^0(EG\times M)\otimes(J\C G)^n$ for all $n$, hence extends to an homomorphism of completed algebras
\begin{equation}
\rho_*\ :\ \Th\A \to \Hh\ .
\end{equation}
We recall that the (completed) space of non-commutative differential forms $\Omh\Th\A=\prod_{n=0}^{\infty}\Omega^n\Th\A$ endowed with the total differential $(b+B)$ computes the periodic cyclic homology of $\A$. \\

The next step is a slight generalization of the Cuntz-Quillen $X$-complex \cite{CQ95} to the DG algebra setting. Indeed the de Rham differential $d$ on $\Hh$ extends in a unique way to the $\Hh$-bimodule of universal 1-forms $\Omega^1\Hh$ by
\[
d(\hat{h}_0\dd\hat{h}_1) = (d\hat{h}_0)\dd\hat{h}_1 + (-1)^{|\hat{h}_0|+1}\hat{h}_0\dd(d\hat{h}_1) \ ,\qquad \forall \ \hat{h}_0\dd\hat{h}_1\in \Omega^1\Hh\ ,
\]
where $|\hat{h}_0|$ denotes the degree of $\hat{h}_0$. Then $d$ is a differential of odd degree on $\Omega^1\Hh$ endowed with its total grading, compatible with the bimodule structure, and commutes in the graded sense with the universal differential $\dd$. We define the $X$-complex of the DG algebra $(\Hh,d)$ as the $\Z_2$-graded supercomplex
\begin{equation}
X(\Hh,d)\ :\ \Hh \rightleftarrows \Omega^1\Hh_\natural\ ,
\end{equation}
where $\Omega^1\Hh_\natural = \Omega^1\Hh/[\Hh,\Omega^1\Hh]$ is the quotient of the bimodule of universal 1-forms by its subspace of graded commutators. We denote by $\natural \hat{h}_0\dd\hat{h}_1$ the class of $\hat{h}_0\dd\hat{h}_1$. The map $\natural\dd: \Hh\to \Omega^1\Hh_\natural$ is simply $\hat{h}\mapsto \natural\dd\hat{h}$, while $\overline{b}: \Omega^1\Hh_\natural\to \Hh$ descends from the graded Hochschild boundary operator $\hat{h}_0\dd\hat{h}_1\mapsto (-1)^{|\hat{h}_0|}[\hat{h}_0,\hat{h}_1]$. One has $\natural\dd\circ\overline{b} = 0$, $\overline{b}\circ \natural\dd =0$, and the odd differential $\natural\dd\oplus\overline{b}$ commutes in the graded sense with $d$, so that $X(\Hh,d)$ endowed with the total differential $(\natural\dd\oplus\overline{b}) + d$ is a $\Z_2$-graded complex. The proof of the following lemma is a straightforward computation.
\begin{lem}\label{lchi}
The linear map of even degree $\chi(\rho_*,d)$ from $\Omh\Th\A$ to $X(\Hh,d)$ given by
\begin{eqnarray}
\lefteqn{\chi(\rho_*,d)(\ah_0 \dd\ah_1\ldots \dd\ah_n)=} \\
&& \frac{1}{(n+1)!} \sum_{i=0}^n (-1)^{i(n-i)} d\rho_*(\ah_{i+1}) \ldots d\rho_*(\ah_n) \, \rho_*(\ah_0) \, d\rho_*(\ah_1) \ldots d\rho_*(\ah_i) \nonumber \\
&& + \frac{1}{n!} \sum_{i=1}^n \natural\big(\rho_*(\ah_0) \, d\rho_*(\ah_1) \ldots \dd \rho_*(\ah_i) \ldots d\rho_*(\ah_n) \big) \nonumber
\end{eqnarray}
for all $\ah_i\in\Th\A$, is a cocycle in the $\Hom$-complex $\Hom(\Omh\Th\A, X(\Hh,d))$. 
\end{lem}
Notice that a differential form $d\rho_*(\ah)\in\Hh$ has always degree 0 in the direction $EG$ and degree 1 in the direction $M$, so that $\chi(\rho_*,d)$ vanishes on $\Omega^n\Th\A$ whenever $n>\dim M$ and thus extends to the direct product $\Omh\Th\A=\prod_{n=0}^{\infty}\Omega^n\Th\A$. This would not be the case if the image of $\A$ in $\G$ consisted in non-constant functions in the direction $EG$.\\

The last step associates a cocycle $\lambda_\omega'\in \Hom(X(\Hh,d),\C)$ to any closed $G$-invariant differential form $\omega\in \Omega(EG\times M)$. To that purpose we define the $X$-complex \emph{localized at units} as the vector space
\begin{equation}
X(\Hh,d)_{[EG\times M]} = \bigoplus_{k\geq 0} \varprojlim_n \big( \Omega_p^k(EG\times M)\otimes \Omega^n\C G_{[1]}\big) \ ,
\end{equation}
where $\Omega^n\C G_{[1]}$ is the space of universal $n$-forms localized at the unit $1\in G$, with the additional relation that $dU_g=0$ for $g=1$ :
\[
U_{g_0}dU_{g_1}\ldots dU_{g_n} \in \Omega^n\C G_{[1]} \Leftrightarrow g_0g_1\ldots g_n =1\ , 
\]
\[
dU_{g_1}\ldots dU_{g_n} \in \Omega^n\C G_{[1]} \Leftrightarrow g_1\ldots g_n =1\ . 
\]
$X(\Hh,d)_{[EG\times M]}$ is a quotient of $X(\Hh,d)$. Indeed the projection $c:\H\to\Omega_p(EG\times M)\otimes \Omega^+\C G_{[1]} $ is unambiguously defined by (see \cite{CQ95} for an explanation of the normalization factor of $(-1)^nn!$)
\begin{multline}
c(\alpha\otimes U_{g_0}\otimes (U_{g_1g_2}-U_{g_1}U_{g_2})\otimes\ldots\otimes (U_{g_{2n-1}g_{2n}}-U_{g_{2n-1}}U_{g_{2n}})) \nonumber\\
= (-1)^nn!\, \alpha\otimes U_{g_0}dU_{g_1}dU_{g_2}\ldots dU_{g_{2n-1}}dU_{g_{2n}}\nonumber
\end{multline}
whenever $g_0g_1\ldots g_{2n} =1$, 
\begin{multline}
c(\alpha\otimes (U_{g_1g_2}-U_{g_1}U_{g_2})\otimes\ldots\otimes (U_{g_{2n-1}g_{2n}}-U_{g_{2n-1}}U_{g_{2n}})) \nonumber\\
= (-1)^nn!\, \alpha\otimes dU_{g_1}dU_{g_2}\ldots dU_{g_{2n-1}}dU_{g_{2n}}\nonumber
\end{multline}
whenever $g_1\ldots g_{2n} =1$, and $c$ vanishes on all other tensors. In the same manner $c:\Omega^1\H_\natural \to\Omega_p(EG\times M)\otimes \Omega^-\C G_{[1]}$ is uniquely specified by
\begin{multline}
c(\natural (\alpha\otimes U_{g_0}\otimes (U_{g_1g_2}-U_{g_1}U_{g_2})\otimes\ldots\otimes (U_{g_{2n-1}g_{2n}}-U_{g_{2n-1}}U_{g_{2n}}))\dd (\beta\otimes U_{g_{2n+1}})) \nonumber\\
= (-1)^{n+|\beta|} n!\, \alpha\wedge U_{g_0\ldots g_{2n}}(\beta) \otimes U_{g_0}dU_{g_1}dU_{g_2}\ldots dU_{g_{2n-1}}dU_{g_{2n}}dU_{g_{2n+1}}\nonumber
\end{multline}
whenever $g_0g_1\ldots g_{2n+1} =1$, 
\begin{multline}
c(\natural (\alpha\otimes (U_{g_1g_2}-U_{g_1}U_{g_2})\otimes\ldots\otimes (U_{g_{2n-1}g_{2n}}-U_{g_{2n-1}}U_{g_{2n}}))\dd (\beta\otimes U_{g_{2n+1}})) \nonumber\\
= (-1)^{n+|\beta|} n!\, \alpha\wedge U_{g_1\ldots g_{2n}}(\beta) \otimes dU_{g_1}dU_{g_2}\ldots dU_{g_{2n-1}}dU_{g_{2n}}dU_{g_{2n+1}}\nonumber
\end{multline}
whenever $g_1\ldots g_{2n+1} =1$, and $c$ vanishes in all other cases. $|\beta|$ is the degree of the differential form $\beta$. By construction $c$ is compatible with the filtration by the powers of the ideal $J\C G$, hence extends to a well-defined map $X(\Hh,d)\to X(\Hh,d)_{[EG\times M]}$. Moreover the boundary operators $\natural\dd$, $\overline{b}$ and $d$ descend to boundary operators on the quotient. Thus $X(\Hh,d)_{[EG\times M]}$ is a quotient of the complex $X(\Hh,d)$. For any $G$-invariant differential form $\omega\in \Omega(EG\times M)$ we define a cochain $\lambda_\omega \in\Hom(X(\Hh,d)_{[EG\times M]} , \C)$ by 
\begin{equation}
\lambda_\omega (\alpha\otimes U_{g_0}dU_{g_1}\ldots dU_{g_n}) = \frac{1}{(2\pi \i)^m}
\int_{\widetilde{\Delta}(g_1,\ldots, g_n)\times M} \alpha \wedge \omega \ ,\qquad m=\dim\, M \ ,\label{lamb}
\end{equation}
where $\widetilde{\Delta}(g_1,\ldots, g_n) \subset EG$ denotes the $n$-simplex with vertices $g_1\ldots g_n$, $g_2\ldots g_n$, $\ldots$, $g_n$, $1$. The normalization factor involving powers of $2\pi\i$ is inserted for ensuring compatibility with index pairings. Remark that the integral above is well-defined, because $\alpha$ has compact $M$-support at any point of $EG$. One can be worried by the fact that by definition of the completion $\Hh$, the degree $k$ of the differential form $\alpha$ is fixed while $n$ can be arbitrarily large. But this causes no trouble since the r.h.s. of (\ref{lamb}) vanishes for large $n$. A direct computation involving Stokes theorem yields

\begin{lem}\label{llambda}
The map sending any $G$-invariant differential form $\omega\in \Omega(EG\times M)$ to the cochain $\lambda_\omega' = \lambda_\omega\circ c \in \Hom(X(\Hh,d),\C)$ is a morphism of complexes.
\end{lem}

Recall that the $\Hom$-complex $\Hom(\Omh\Th\A,\C)$ computes the periodic cyclic cohomology of $\A$. Collecting Lemmas \ref{lchi} and \ref{llambda} one thus gets

\begin{prop}
Let $G$ be a discrete group acting by orientation-preserving diffeomorphisms on a smooth oriented manifold $M$, and $\A = \cinfc(M)\rtimes G$. The map sending any $G$-invariant differential form $\omega\in \Omega(EG\times M)$ to the cochain $\lambda_\omega'\circ \chi(\rho_*,d) \in \Hom(\Omh\Th\A,\C)$ is a morphism of complexes. We denote by
\begin{equation}
\Phi : H^\bullet(EG\times_G M) \to \HP^\bullet(\A)
\end{equation}
the corresponding map in cohomology.
\end{prop}

\section{Algebraic JLO formula}\label{sJLO}

In this section we collect all the preceding results and construct  cyclic cocycles by means of JLO-type formulas. Let $M$ be a foliated manifold of dimension $n=v+h$, $v$ being the dimension of the leaves, and $G\subset\Diff(M)$ a discrete group of foliated diffeomorphisms. Let $EG$ be the universal $G$-bundle over the classifying space $BG$. We take $EG$ as the geometric realization of the simplicial set $N\overline{G}_{\bullet}$ and consider the induced action of $G$ on the DG algebra of smooth differential forms $\Omega(EG)$. We shall denote by $d_H$ the differential over $EG$, reserving the classical notation $d$ for the total differential over the products of spaces $EG\times M$ or $EG \times T^*M$. The algebraic crossed product
\begin{equation}
\G = \Omega(EG)\rtimes G
\end{equation}
is a particular case of the DG algebra constructed in Section \ref{sequiv}, the manifold $M$ being replaced by a point. Hence considering twisted tensor products with the tensor algebra of $\C G$, one gets the DG algebra extension of $\G$ and its completion:
\begin{equation}
\H = \Omega(EG)\otimes T\C G \ ,\qquad \Hh = \bigoplus_{k\geq 0}\varprojlim_m \big(\Omega^k(EG)\otimes T\C G/(J\C G)^m\big) \ .
\end{equation}
Now let $\S(M)=\L(M)[[\eps]]$ be the $G$-algebra of formal power series constructed in Section \ref{sbimod}, together with the subalgebras $\D(M)$, $\D_c(M)$, and the corresponding $\D(M)$-bimodule of trace-class operators $\T(M)$. The space $\Omega(EG,\S(M))$ of smooth differential forms on $EG$ with values in $\S(M)$ is a DG algebra, for the pointwise product of differential forms and de Rham differential $d_H$ over $EG$. We endow $\Omega(EG,\S(M))$ with the $G$-action which combines the actions of $G$ on $EG$ and $\S(M)$ respectively. The crossed product 
\begin{equation}
\U = \Omega(EG,\S(M))\rtimes G
\end{equation}
is therefore a DG algebra, with differential $d_H(\alpha\otimes U_g) = d_H\alpha\otimes U_g$ for all $\alpha\in \Omega(EG,\S(M))$ and $g\in G$. Considering twisted tensor products with the tensor algebra of $\C G$, a DG algebra extension of $\U$ is defined as above:
\begin{equation}
\V = \Omega(EG,\S(M))\otimes T\C G \ .
\end{equation}
We will \emph{not} take the adic completion of $\V$ with respect to the kernel of the projection homomorphism $\V\to \U$ because this does not lead to the right object. Instead, we only build an homomorphism from the tensor algebra $T\A$ over $\A=\cinfc(S^*_HM)\rtimes G$ to $\V$. As in Section \ref{sindex} choose a linear splitting 
\[ \nu: \A \to \Psi^0_{H,c}(M) \rtimes G \]
of the leading symbol homomorphism. Then observe that the trivial line bundle $M\times\C$ over $M$ can be identified with the zero-degree part of the exterior bundle $E=\Lambda^\bullet( T^* M\otimes\C)$. This allows to view the scalar symbols $\S_{H,c}(M)$ as a subalgebra of $\S_{H,c}(M,E)$. One thus has a sequence of $G$-equivariant homomorphisms
\[ \Psi^0_{H,c}(M) \to \S^0_{H,c}(M) \hookrightarrow \S^0_{H,c}(M,E) \stackrel{L}{\longrightarrow} \D_c(M) \hookrightarrow \Omega^0(EG,\S(M)) \]
where $L$ is the representation of symbols by left multiplication, and $\D_c(M)$ is viewed as an algebra of constant functions on $EG$. Composing $\nu$ with the above sequence leads to a linear map $\sigma:\A \to \U$ which, by the universal property of the tensor algebra, extends to an homomorphism of algebras 
\begin{equation}
\sigma_*\ :\ T \A \to \V\ .
\end{equation}
Of course the latter depends on the original choice of linear splitting $\nu$, but two different splittings lead to homotopic homomorphisms in the sense of Cuntz and Quillen \cite{CQ95}.\\

We now discuss superconnections \cite{Qui1988}. Remark that the DG algebra $\Omega(EG,\D(M))$ acts by left and right multipliers on $\U$ and $\V$. Let $D_0\in \D(M)$ be a generalized Dirac operator as defined in Section \ref{sdirac} and consider the element $D\in \Omega^0(EG,\D(M))$ as the function on the classifying space with values in the set of Dirac operators, given by
\[
D(g_0,\ldots,g_r)(s_0,\ldots,s_r) = \sum_{i=0}^r s_i \, \Ad_{g_i}^{-1}(D_0)
\]
for all $(g_0,\ldots,g_r)\in N\overline{G}_r$ and $(s_0,\ldots,s_r)\in \Delta_r$. The action of an element $g\in G$ on $EG$ carries $(g_0,\ldots,g_r)$ to $(g_0g,\ldots,g_rg)$, so the function $D$ is $G$-equivariant by construction. The superconnection 
\begin{equation}
\DD = d_H + D\ ,
\end{equation}
acting on $\V$ by graded commutators, is a graded derivation. Its curvature is the inhomogeneous differential form 
\[
\DD^2 = D^2 + d_HD  \ \in \Omega^0(EG,\D(M))\oplus \Omega^1(EG,\D(M))
\]
where $d_HD(g_0,\ldots,g_r)(s_0,\ldots,s_r) = \sum_{i=0}^r ds_i \, \Ad_{g_i}^{-1}(D_0)$. 
Choose a positive Heisenberg-elliptic symbol $q\in \S_H^1(M)$ of order one on $M$. Typically, we take $q$ as the symbol of the operator $\Delta_H^{1/4}$ associated to a hypoelliptic sub-Laplacian. Extend it to a Heisenberg-elliptic symbol $\tilde{q}\in \S_H^1(M,E)$, requiring that the leading symbol of $\tilde{q}$ remains of scalar type. Using the left representation of symbols $L:\S_H(M)\to \D(M)$, one gets a constant function $\tilde{q}_L \in \Omega^0(EG, \D(M))$ over $EG$. Let $\kappa$ be an ``infinitesimal'' odd parameter: $\kappa^2=0$. The new superconnection
\begin{equation}
\nabla = \DD + \kappa\ln\tilde{q}_L\ ,
\end{equation}
acting on the DG algebra $\V[\kappa]= \V\oplus\kappa\V$ by graded commutators, is a graded derivation. Its curvature is the inhomogeneous differential form 
\[
\nabla^2 = \DD^2 + \kappa [\ln\tilde{q}_L,\DD] \ \in \Omega^0(EG,\D(M))[\kappa]\oplus \Omega^1(EG,\D(M))[\kappa]\ .
\]

The homomorphism $\sigma_*$ and the superconnection $\nabla$ are the main ingredients of a JLO-type formula for a cocycle of odd degree in $\Hom(\Omh\Th\A, X(\Hh,d_H)_{[EG]})$, where $X(\Hh,d_H)_{[EG]}$ is the $X$-complex \emph{localized at units} defined in Section \ref{sequiv}. We first extend the graded trace  $\Tr_s:\T(M)\to\C$ of Section \ref{scano} to an $X$-complex map in an obvious way :
\[ \Tr_s : X(\V) \to X(\H,d_H) \]

\begin{prop}
The linear map of even degree $\chi^{\Tr_s}(\sigma_*,\nabla) : \Omega T\A \to X(\H,d_H)[\kappa]$ defined on any $r$-form $\ah_0\dd\ah_1\ldots\dd\ah_r \in \Omega^r T\A$ by
\begin{eqnarray}
\lefteqn{\chi^{\Tr_s}(\sigma_*,\nabla)(\ah_0\dd\ah_1\ldots\dd\ah_r)=}\\
&&\hspace{-1cm}\sum_{i=0}^r (-)^{i(r-i)} \int_{\Delta_{r+1}} \Tr_s\big(e^{-t_{i+1}\nabla^2} [\nabla,\sigma_{i+1}] \ldots  e^{-t_{r+1}\nabla^2} \sigma_0 \,e^{-t_0\nabla^2} [\nabla,\sigma_1] \ldots e^{-t_i\nabla^2}\big) dt \nonumber\\
&+&  \sum_{i=1}^r \int_{\Delta_r} \Tr_s\big(\natural \sigma_0\, e^{-t_0\nabla^2} [\nabla,\sigma_1] \ldots  e^{-t_{i-1}\nabla^2} \dd \sigma_i \,e^{-t_i\nabla^2} \ldots [\nabla,\sigma_r] e^{-t_r\nabla^2}\big)dt \nonumber
\end{eqnarray}
where $\sigma_i = \sigma_*(\ah_i)\in \V$ for all $i$, extends to a cochain in the complex $\Hom(\Omh\Th\A, X(\Hh,d_H))[\kappa]$. Moreover, its composition with the projection map onto the $X$-complex localized at units $c:X(\Hh,d_H) \to X(\Hh,d_H)_{[EG]}$ yields a \emph{cocycle} 
\[ c\circ\chi^{\Tr_s}(\sigma_*,\nabla) \in \Hom(\Omh\Th\A, X(\Hh,d_H)_{[EG]})[\kappa] \ .\]
\end{prop}
The above expression is well-defined, because it involves a Duhamel-type expansion of the heat operator $\exp(-\nabla^2)$ which belongs to the domain of the trace. The extension of the cocycle $\chi^{\Tr_s}(\sigma_*,\nabla)$ to completions follows from a straightforward adaptation of the construction of the bivariant Chern character given in \cite{Per2004}. In our case, things are even easier since we do not have to work in entire cyclic cohomology. That $c\circ \chi^{\Tr_s}(\sigma_*,\nabla)$ is a cocycle follows from well-known algebraic manipulations \cite{JLO}, which in the case at hand crucially depend on the formal identities $\dd D=0$ and $\dd\ln\tilde{q}_L=0$ holding only in the localized complex. Since $\kappa^2 = 0$, the cochain $\chi^{\Tr_s}(\sigma_*,\nabla)$ is actually a polynomial of degree one with respect to $\kappa$. Define
\begin{equation}
\chi^{\Tr_s}(\sigma_*,\DD,\ln\tilde{q}_L) = \frac{\partial}{\partial \kappa} \chi^{\Tr_s}(\sigma_*,\DD+\kappa\ln\tilde{q}_L)\ .
\end{equation}
The latter does no longer depend on $\kappa$ and is a cochain of odd degree in $\Hom(\Omh\Th\A,X(\Hh,d_H))$. Likewise, the composite 
\[ c\circ \chi^{\Tr_s}(\sigma_*,\DD,\ln\tilde{q}_L) \in \Hom(\Omh\Th\A,X(\Hh,d_H)_{[EG]}) \]
is a \emph{cocycle} of odd degree. By a classical homotopy argument, its cohomology class does not depend on any choice regarding the linear map $\sigma:\A\to\U$, the superconnection $\DD$, and the elliptic symbol $\tilde{q}$. Thus, \emph{this class is completely canonical and only depends on the data $(M,G)$}. The following proposition identifies its composition with the class of the chain map $\lambda_1 \in\Hom(X(\Hh,d_H)_{[EG]},\C)$ of Lemma \ref{llambda}, for $\omega = 1$.

\begin{prop} \label{pDeRham}
Let $D_0\in \D(M)$ be a de Rham-Dirac operator, $D\in \Omega^0(EG,\D(M))$ the associated $G$-equivariant family of Dirac operators on the universal bundle, and $\DD = d_H + D $ the corresponding superconnection. Then $\lambda'_1\circ \chi^{\Tr_s}(\sigma_*,\DD,\ln \tilde{q}_L)$ is the cocycle of Proposition \ref{pcocycle}. 
\end{prop}
\begin{pr}
We work in a local foliated chart $(x,p)$ over $T^*M$. First, we must observe that $D$ is still a de Rham - Dirac type operator, essentially because $d_R = \i p_i \psi^{i}$ is $G$-invariant. Thus, $D$ is of the form 
\begin{align*}
D = -\i \eps d_R + \overline{\psi}_{i R} \left(\partial_{p_{i}} + \sum_{\vert \alpha \vert \geq 2} (r^{i}_{\alpha})_L  \partial^{\alpha}_p \right)  
\end{align*}
The $r^{i}_{\alpha}$ are scalar functions on $EG\times M$. With $d_H$ the exterior differential on $EG$, one has
\[ d_HD = \sum_{\vert \alpha \vert \geq 2} (d_Hr^{i}_{\alpha})_L \overline{\psi}_{i R}  \partial^{\alpha}_p\]
Moreover, $D^{2}$ reads 
\begin{multline} \label{laplacian bis}
-D^{2} = \Delta  + \eps  (p_{iL} \partial_{p_{i}} + (\psi^{i} \overline{\psi}_i)_R) \\
+ \eps \left(\sum_{\vert \alpha \vert \geq 2} ((a_{\alpha}^{i})_L \partial_{x^{i}} + (a_{\alpha}^{i} p_i)_L) \partial^{\alpha}_{p} + \sum_{\vert \alpha \vert \geq 1} (b_{\alpha j}^{i})_L (\psi^{j} \overline{\psi}_i)_R \partial^{\alpha}_{p} \right)
\end{multline}
where $\Delta = \i \eps \partial_{x^{i}} \partial_{p_{i}}$, and the coefficients in the sums over $\vert \alpha \vert \geq \ldots$ are also scalar functions on $EG\times M$. Recall also that
\begin{align*}
\DD = d_H + D, && \DD^{2} = d_HD + D^{2}, && \nabla = \DD + \kappa \ln \tilde{q}_L, && \nabla^{2} = \DD^{2} + \kappa [\ln \tilde{q}_L, \DD] 
\end{align*}
We know that $d_HD$ is proportional to $\psio_R$. The symbol $q$ being constant in the direction $EG$, one has $d_H\ln \tilde{q}_L=0$. Moreover $\ln \tilde{q}_L$ commutes with $p_R$, hence the commutator $[\ln \tilde{q}_L, \DD]$ is also proportional to $\psio_R$. Finally the elements $\sigma=\sigma_*(\ah)\in \Omega^0(EG,\D_c(M))\otimes T\C G$ are constant in the direction $EG$ for all $\ah\in T\A$, hence we have 
\begin{equation*}
[\nabla, \sigma] = [\overline{\psi}_{i R} (\partial_{p_{i}} + \ldots), \sigma] + \kappa [\ln \tilde{q}_L , \sigma] \ .
\end{equation*} 
Now observe that the graded trace $\Tr_s$ selects the term proportional to $(\psi^{1} \overline{\psi}_1 \ldots \psi^{n} \overline{\psi}_n)_R$. The generalized Laplacian $D^{2}$ already brings terms proportional to $1$ or $(\psi \overline{\psi})_R$ in the right sector. Thus the terms proportional to $\overline{\psi}_{R}$ in $d_HD$, $[\ln \tilde{q}_L, \DD]$ and $[\nabla, \sigma]$ break the balance between the $\psi_R$ and the $\overline{\psi}_R$ and must give a zero contribution to the cochain $\chi^{\Tr_s}(\sigma_*,\nabla)$. Hence we can consider that
\[ \nabla^{2} \simeq D^{2} \ ,\qquad [\nabla, \sigma] \simeq  \kappa [\ln \tilde{q}_L , \sigma]\ .\]
A first consequence, taking into account $\kappa^{2} = 0$, is that the cochain $\chi^{\Tr_s}(\sigma_*,\DD,\ln\tilde{q}_L)$ should contain exactly one commutator $[\nabla, \sigma]$. Thus the only non-zero components of this cochain are
\begin{align*}
\chi^{\Tr_s}(\sigma_*,\DD,\ln\tilde{q}_L)(\ah_0 \textbf{d}\ah_1) =& \int_{\Delta_2} \Tr_s \left( e^{-t_1 D^{2}} [\ln \tilde{q}_L,\sigma_1] e^{-t_2 D^{2}} \sigma_0 e^{-t_0 D^{2}} \right) dt \\ 
&+ \int_{\Delta_2} \Tr_s \left( e^{-t_2 D^{2}} \sigma_0 e^{-t_1 D^{2}} [\ln \tilde{q}_L, \sigma_1] e^{-t_0 D^{2}} \right) dt 
\end{align*}
\begin{align*}
\chi^{\Tr_s}(\sigma_*,\DD,\ln\tilde{q}_L)(\ah_0 \textbf{d}\ah_1\dd \ah_2) =& \int_{\Delta_2} \Tr_s \left(\natural \sigma_0 e^{-t_0 D^{2}} [\ln \tilde{q}_L, \sigma_1] e^{-t_1 D^{2}} \dd\sigma_2 e^{-t_2 D^{2}} \right) dt \\ 
&+ \int_{\Delta_2} \Tr_s \left(\natural \sigma_0 e^{-t_0 D^{2}} \dd\sigma_1 e^{-t_1 D^{2}} [\ln \tilde{q}_L, \sigma_2] e^{-t_2 D^{2}} \right) dt 
\end{align*}
A second consequence is that the images of these quantities under the projection $c:X(\Hh,d_H)\to X(\Hh,d_H)_{[EG]}$ belong to the subspace $\Omega^0(EG)\otimes\Omega\C G_{[1]}$ of the localized $X$-complex, in other words they are scalar functions over $EG$. Thus, their evaluation on the cocycle $\lambda_1$ drops all the components in $\Omega^0(EG)\otimes\Omega^k\C G$ for $k\geq 1$, and the remaining components in $\Omega^0(EG)\otimes\Omega^0\C G$ are simply localized at the basepoint $(1,1)\in EG\times G$. In particular 
\[
\lambda_1'\circ\chi^{\Tr_s}(\sigma_*,\DD,\ln\tilde{q}_L)(\ah_0 \textbf{d}\ah_1\dd \ah_2) = 0
\]
and \emph{$\lambda_1'\circ\Tr_s$ behaves like a graded trace} in the only remaining term (the subscript $_{[1,1]}$ denotes localization at the basepoint) : 
\begin{align*} 
\lambda_1'\circ\chi^{\Tr_s}(\sigma_{*}, \DD,\ln\tilde{q}_L)(\ah_0 \textbf{d}\ah_1) =& \int_{\Delta_2} \Tr_s \left( e^{-t_1 D^{2}} [\ln\tilde{q}_L, \sigma_1] e^{-t_2 D^{2}} \sigma_0 e^{-t_0 D^{2}} \right)_{[1,1]} dt \\ 
& + \int_{\Delta_2} \Tr_s \left( e^{-t_2 D^{2}} \sigma_0 e^{-t_1 D^{2}} [\ln\tilde{q}_L, \sigma_1] e^{-t_0 D^{2}} \right)_{[1,1]} dt \\
=& \int_0^{1} \Tr_s \left( \sigma_0 e^{-t D^{2}} [\ln\tilde{q}_L, \sigma_1] e^{-(1-t) D^{2}} \right)_{[1,1]} dt \\
\end{align*}
The integrand $\Tr_s (\sigma_0 e^{-t D^{2}} [\ln\tilde{q}_L, \sigma_1] e^{-(1-t) D^{2}} )_{[1,1]}$ does not depend on $t \in [0,1]$. Indeed 
\begin{align*}
\!\!\!\!\!\!\! \dfrac{d}{dt} \Tr_s \left( \sigma_0 e^{-t D^{2}} [\ln\tilde{q}_L, \sigma_1] e^{-(1-t) D^{2}} \right)_{[1,1]} & = - \Tr_s \left(\sigma_0 e^{-t D^{2}} [D^2,[\ln\tilde{q}_L, \sigma_1]] e^{-(1-t) D^{2}} \right)_{[1,1]} \\
& = \Tr_s \left( [D,\sigma_0] e^{-t D^{2}} [D,[\ln\tilde{q}_L, \sigma_1]] e^{-(1-t) D^{2}} \right)_{[1,1]}
\end{align*}
The last equality comes from $[D^{2}, X] = D [D, X] + [D, X] D $ and the graded trace property. The above quantity vanishes,  because the commutators with $D$ only bring terms proportional to $\psio_R$. Therefore, the integrand may be replaced by its value at $t=0$, and we are left with 
\[ \lambda_1'\circ\chi^{\Tr_s}(\sigma_{*}, \DD,\ln\tilde{q}_L)(\ah_0 \textbf{d}\ah_1) = \Tr_s \left( \sigma_0 [\ln \tilde{q}_L, \sigma_1] e^{-D^{2}} \right)_{[1,1]} \]
Seeing $-D^{2}$ as a perturbation of the flat Laplacian $\Delta + u$, and using a Duhamel expansion, one gets 
\begin{align*}
\Tr_s \left( \sigma_0 [\ln \tilde{q}_L, \sigma_1] e^{-D^{2}} \right) = \sum_{k \geq 0} \int_{\Delta_k} \Tr_s \left( \sigma_0 [\ln \tilde{q}_L, \sigma_1]\sigma_{\Delta}^{t_0}(u) \ldots \sigma_{\Delta}^{t_0 + \ldots + t_{k-1}}(u) \exp(\Delta) \right) dt 
\end{align*}
Then, we want to move the operators $\partial_x$ and $\partial_p$ to the right in each term of the sum above, and look at when we have an exact balance in their powers. Otherwise, it will vanish under the graded trace $\Tr_s$ by definition. A $\partial_p$ can be absorbed with a $p_L$ by commutation, and $\partial_x$ may appear in $\sigma_\Delta^{t}(p_L) = p_L + t[\Delta, p_L] )=  p_L + it\eps \partial_x$. With this elements at hand, we conclude that the presence of the sums over $\vert \alpha \vert$ in (\ref{laplacian bis}) prevent an exact balance between $\partial_x$ and $\partial_p$. So, we can neglect these parts in $D^{2}$ and get      
\begin{align*}
\!\!\!\!\!\!\!\! \Tr_s \left( \sigma_0 [\ln \tilde{q}_L, \sigma_1] e^{-D^{2}} \right)_{[1,1]} & = \Tr_s \left( \sigma_0 [\ln \tilde{q}_L, \sigma_1] \exp(\Delta + \eps p_L \cdot \partial_p + \eps (\psi^{i} \psio_i)_R) \right)_{[1,1]} \\
& = \Tr_s\left( \sigma_0 [\ln \tilde{q}_L, \sigma_1] \eps^{n} (\psi^{1} \psio_1 \ldots \psi^{n} \psio_n)_R exp(\Delta + \eps p_L \cdot \partial_p) \right)_{[1,1]} 
\end{align*}
where in the second equality, we split the exponential. Let $\nu: \A\to \Psi^0_{H,c}(M)\rtimes G$ be the linear splitting that we used in the construction of $\sigma$, and extend it to an algebra homomorphism $\nu_*: T\A \to \Psi^0_{H,c}(M)\rtimes G$. Then setting $\nu_i=\nu_*(\ah_i)$ we get
\begin{align*}
\!\!\!\!\!\!\!\! \Tr_s \left( \sigma_0 [\ln \tilde{q}_L, \sigma_1] e^{-D^{2}} \right)_{[1,1]} & = \barint (\nu_0 [\ln q, \nu_1])_{[1]} \langle \langle \eps^n (\psi^{1} \psio_1 \ldots \psi^{n} \psio_n)_R exp(\Delta + \eps p_L \cdot \partial_p) \rangle \rangle[n] \\
& = \barint (\nu_0 [\ln q, \nu_1])_{[1]} \langle exp(\Delta + \eps p_L \cdot \partial_p) \rangle [0]
\end{align*}
where $_{[1]}$ denotes localization at the unit $1\in G$. Example \ref{todd matrix} applied to the matrix $\Omega=\eps\Id$ yields $\langle exp(\Delta + \eps p_L \cdot \partial_p) \rangle [0]=1$, hence
\[ \Tr_s \left( \sigma_0 [\ln \tilde{q}_L, \sigma_1] e^{-D^{2}} \right)_{[1]} =  \barint (\nu_0 [\ln q, \nu_1])_{[1]} \]
is the equivariant Radul cocycle of Proposition \ref{pcocycle}. $\hfill{\square}$
\end{pr}

Finally consider the diagonal action of $G$ on the product $EG\times S^*_HM$, and its induced action on the space of differential forms $\Omega(EG\times S^*_HM)$ with total de Rham differential $d$. As in Section \ref{sequiv} we form the DG algebra
\begin{equation}
\X = \Omega_p(EG\times S^*_HM)\rtimes G\ ,
\end{equation}
and its DG extension 
\begin{equation}
\Y = \Omega_p(EG\times S^*_HM)\otimes T\C G\ ,\qquad \Yh = \bigoplus_{k\geq 0}\varprojlim_n \big(\Omega_p^k(EG\times S^*_HM)\otimes T\C G/(J\C G)^n\big) \ .
\end{equation}
Viewing $\cinfc(S^*_HM)\subset \Omega_p^0(EG\times S^*_HM)$ as the subalgebra of constant functions in the direction $EG$ leads to an homomorphism $\rho:\A\to\X$ which extends to
\begin{equation}
\rho_* \ :\ \Th\A \to \Yh\ .
\end{equation}
Let $\chi(\rho_*,d) \in \Hom(\Omh\Th\A,X(\Yh,d))$ be the cocycle of Lemma \ref{lchi}. The integration of differential forms along the fibers of the projection $EG\times S^*_HM \to EG$ yields a morphism of complexes 
\begin{equation}
\int_{S^*_HM} \ :\ X(\Yh,d) \to X(\Hh,d_H)\ .
\end{equation}
For the Proposition below we choose a positive Heisenberg-elliptic symbol $q_0\in \S_H^1(M)$ of order one on $M$, and extend it to a Heisenberg-elliptic symbol $\tilde{q}_0\in \S_H^1(M,E)$, requiring that the leading symbol of $\tilde{q}_0$ remains of scalar type. Then the function $\tilde{q}\in \Omega^0(EG,\S_H^1(M,E))$ defined on the universal $G$-bundle by
\[
\tilde{q}(g_0,\ldots,g_r)(s_0,\ldots,s_r) = \sum_{i=0}^r s_i \, \Ad_{g_i}^{-1}(\tilde{q}_0)
\]
is $G$-equivariant. The corresponding function $\tilde{q}_L\in \Omega^0(EG,\D(M))$ is therefore also $G$-equivariant. 

\begin{prop} \label{pAffine}
Let $D_0\in \D(M)$ be a generalized Dirac operator affiliated to a torsion-free affine connection $\Gamma_0$ on $M$. Let $D\in \Omega^0(EG,\D(M))$ be the associated $G$-equivariant family of Dirac operators over $EG$, and $\DD = d_H + D $ the corresponding superconnection. Then in $\Hom(\Omh\Th\A,X(\Hh,d_H))$ holds the equality of cochains
\begin{equation}
\chi^{\Tr_s}(\sigma_*,\DD,\ln \tilde{q}_L) = \frac{1}{(2\pi \i)^n} \int_{S^*_HM} \chi(\rho_*,d) \wedge \Td(R) \ ,\qquad n=\dim\, M\ ,
\end{equation}
where $R$ is the equivariant curvature two-form of $\Gamma$, and $\Td(R)$ is the $G$-invariant closed differential form on $EG\times S^*_HM$ representing the equivariant Todd class of $TM\otimes\C$.
\end{prop}
\begin{pr} First, notice that as a function over $EG$, the operator $D$ takes values in the set of generalized Dirac operators affiliated to affine connections. More precisely, in a foliated local chart, 
\begin{equation*}
D = \i \eps \psi_R^{i} \big(\partial_{x^{i}} + (\Gamma_{ij}^{k})_L(p_{kL} \partial_{p_j} +(\psio_{k} \psi^{j})_L - (\psio_{k} \psi^{j})_R \big) + \psio_{iR} \partial_{p_i} + r, 
\end{equation*}
where the Christoffel symbol $\Gamma_{ij}^{k}$ is a function of the simplicial coordinates $(s_l)$ on $EG$ given by barycentric mean of the affine connection $\Gamma_0$,
\[ \Gamma_{ij}^{k} (g_0,\ldots,g_r)(s_0,\ldots,s_r) = \sum_{l=0}^r s_l \, \Ad_{g_l}^{-1}(\Gamma_0)_{ij}^{k} \ , \]
and $r$ is a remainder of the following form :
\[ r = \i\eps \psi^{i}_R \left(\sum_{\vert \alpha \vert \geq 2} (s_{\alpha i}^{k} p_k)_L \partial_{p}^{\alpha} + \sum_{\vert \alpha \vert \geq 1} (s_{\alpha i j}^{k} \psio_k \psi^{j} + s_{\alpha i})_L \partial_{p}^{\alpha} \right) + \psio_{i R} \sum_{\vert \alpha \vert \geq 2} (r_{\alpha}^{i})_L \partial_{p}^{\alpha}, \]
the coefficients $s^k_{\alpha i}, r^i_\alpha\ldots$ being scalar functions on $EG \times M$. The remainder will in fact not contribute for reasons of order, as we shall explain afterwards. Recall once again that we have 
\begin{align*}
\DD = d_H + D, && \DD^{2} = d_HD + D^{2}, && \nabla = \DD + \kappa \ln \tilde{q}_L, && \nabla^{2} = \DD^{2} + \kappa \delta \DD, 
\end{align*}
where $\delta$ denotes the commutator $[\ln \tilde{q}_L, \ ]$. A straightforward calculation shows that the operator $-\DD^{2}$ is a perturbation of the flat Laplacian $\Delta = \i \eps \partial_{x^{i}} \partial_{p_i}$ (see Proposition \ref{lichner}):
\begin{multline*}
-\DD^{2} =  \Delta + \i \eps (\Gamma_{ij}^{k})_L(\psi^{i} \psio_k)_R \partial_{p_j} \\
 +  \left(-\i\eps d_H(\Gamma_{il}^{k})_L \psi^{i}_R + \dfrac{\eps^2}{2}(R_{lij}^{k})_L (\psi^{i}  \psi^{j})_R \right) \big(p_{kL} \partial_{p_l} + (\psio_k \psi^{l})_L \big) + \ldots     
\end{multline*}
where $R_{lij}^{k}$ denote the components of the curvature tensor of the connection $\Gamma$ along $M$. The dots involve strictly higher powers of $\partial_p$, and will be irrelevant in the calculation, like the remainder $r$ above.  Notice that the coefficient of the third term
\[ \Omega_l^{k}= -\i\eps d_H(\Gamma_{il}^{k})_L \psi^{i}_R + \dfrac{\eps^2}{2}(R_{lij}^{k})_L (\psi^{i}  \psi^{j})_R  \]
 is closely related to the equivariant curvature two-form of the connection $\Gamma$ over $EG\times M$ :
\[ R_l^{k}= d_H\Gamma_{il}^{k} dx^i + \dfrac{1}{2}R_{lij}^{k} dx^{i}\wedge  dx^{j} \ .\]

Finally, recall that
\begin{multline*}
\lefteqn{\chi^{\Tr_s}(\sigma_*,\nabla)(\ah_0\dd\ah_1\ldots\dd\ah_r)=}\\
\sum_{i=0}^r (-)^{i(r-i)} \int_{\Delta_{r+1}} \Tr_s\big(e^{-t_{i+1}\nabla^2} [\nabla,\sigma_{i+1}] \ldots  e^{-t_{r+1}\nabla^2} \sigma_0 \,e^{-t_0\nabla^2} [\nabla,\sigma_1] \ldots e^{-t_i\nabla^2}\big) dt \nonumber\\
+  \sum_{i=1}^r \int_{\Delta_r} \Tr_s\big(\natural \sigma_0\, e^{-t_0\nabla^2} [\nabla,\sigma_1] \ldots  e^{-t_{i-1}\nabla^2} \dd \sigma_i \,e^{-t_i\nabla^2} \ldots [\nabla,\sigma_n] e^{-t_r\nabla^2}\big)dt \nonumber
\end{multline*}
where $\sigma_i = \sigma_*(\ah_i) \in \Omega^0(EG,\D_c(M))\otimes T\C G\subset \V$ for all $a_i\in T\A$. The cochain $\chi^{\Tr_s}(\sigma_*,\DD, \ln \tilde{q}_L)$ is obtained from the latter by extracting the coefficient of $\kappa$. \\

We can now start with the calculations. It suffices to focus on the first sum of the cochain $\chi^{\Tr_s}(\sigma_*, \nabla)$. The other term can be dealt with ad verbatim. \\

The Dirac operator being $G$-equivariant, it commutes to the action of the diffeomorphisms induced by the factor $T\C G$ in the elements $\sigma=\sigma_*(\ah)$. In addition, the leading order part of generalized Dirac operators is $G$-invariant. This $G$-equivariance is crucial all along the proof and allows to perform calculations as if leading terms in the Dirac operator and the Laplacian are $G$-invariant. \\

In particular, the action of $\nabla$ on $\sigma$ via commutators $[\nabla, \sigma]$ reads : 
\[ [\nabla, \sigma] = \i \eps \psi_R^{i} \big(\partial_{x^{i}}\sigma + (\Gamma_{ij}^{k} p_k)_L \partial_{p_{j}} \sigma + \ldots \big) + \psio_{iR} (\partial_{p_i} \sigma + \ldots) + \kappa [\ln \tilde{q}_L, \sigma]  \]
where the term proportional to $\eps \psi_R$ has order $0$, with the dots of order $-1$, and the term proportional to $\psio_{iR}$ has order -1 when $1 \leq i \leq v$, of order -2 when $v+1 \leq i \leq n$. Again, the dots contain higher powers of $\partial_p$, so they have strictly smaller order and will be neglected. \\

Viewing $-\nabla^{2} = \Delta + u$ as a perturbation of the flat Laplacian $\Delta$, a Duhamel expansion in the exponentials leads to studying terms of the form 
\[ \Tr_s \left(\sigma_0 \sigma_{\Delta}^{t_0}(X_1) \ldots \sigma_{\Delta}^{t_0 + \ldots + t_{k-1}}(X_k) \exp \Delta \right)  \]      
where $X_i = u \text{ or } [\nabla, \sigma_j]$. The action of $\sigma_{\Delta}$ does not modify the fact that $X_i$ is of order 0, because 
\[ [\Delta, X_i] = \i \eps \big(\partial_{x^{j}}X_i \cdot \partial_{p_j} + \partial_{p_{j}}X_i \cdot \partial_{x^j} + \partial_{x^{j}} \partial_{p^{j}} X_i \big) \]    
Now, we observe that in the formulas above, $\eps \psi_R$ is always proportional to an operator of order $\leq 0$, $\psio_{iR}$ to an operator of order $\leq -1 \text{ or } -2$ depending on $i$. The graded trace $\Tr_s$ selects the term proportional to $(\psi^{1} \psio_1 \ldots \psi^{n} \psio_n)_R$ which is therefore of order $\leq -(v+2h)$. Since the Wodzicki residue selects terms of order exactly $-(v+2h)$, only the leading terms are involved in these quantities. In particular, the remainder $r$ and the dots in the formulas above will vanish under $\Tr_s$. For a similar reason, the derivatives $\partial_{x}X_i$ that appear in the modular action $\sigma_\Delta$ can be neglected : on the one hand, they do not add up any factor $\psi_R$ or $\psio_R$. On the other hand, they are always proportional to a $\partial_p$ (of order $\leq -1$) . Hence, the factors $\partial_{x^i} X$ decrease the order of terms proportional to $(\psi^{1} \psio_1 \ldots \psi^{n} \psio_n)_R$, and contribute to an overall order $< -(p+2q)$. In other words, when computing the Wodzicki density residue at a point $x_0\in M$ \emph{all functions of $x$ can be considered as constants}. Then, if we choose a coordinate system around $x_0$ such that $\Gamma_{ij}^{k}(x_0) \simeq 0$ over a given point $p_0\in EG$, we have 
\begin{gather*}
-\DD^{2} \simeq \Delta + \Omega_l^{k} \big(p_{kL} \partial_{p_l} + (\psio_k \psi^{l})_L) 
\end{gather*}
and we can ignore the $x$-derivatives of the matrix $\Omega=(\Omega^k_l)$. Another important simplification occurs with the curvature term $\Omega_l^{k} (\psio_k \psi^{l})_L$. In the JLO formula the latter is always multiplied by a factor $\Pi_L$ coming from the $\sigma_i$'s, where $\Pi=\psio_1\psi^1\ldots\psio_n\psi^n$ is the projection operator onto scalar symbols. One has
\[ \Omega_{l}^{k} (\psio_k \psi^{l})_L \Pi_L =  \Omega_{l}^{k} \big((\delta_k^{l} - \psi^{l} \psio_k ) \Pi\big)_L \ ,\qquad \Pi_L \Omega_{l}^{k} (\psio_k \psi^{l})_L  =   \big(\Pi(\delta_k^{l} - \psi^{l} \psio_k) \big)_L\Omega_{l}^{k}  \] 
and the identities $\psio_k \Pi = 0 = \Pi \psi^l$ show that these expressions only depend on the trace $\tr(\Omega)=\Omega_k^k$ of the curvature matrix :
\[ \Omega_{l}^{k} (\psio_k \psi^{l})_L \Pi_L =  \tr(\Omega) \Pi_L \ ,\qquad \Pi_L \Omega_{l}^{k} (\psio_k \psi^{l})_L  =   \Pi_L \tr(\Omega) \ . \] 
Moreover, the commutators $[\tr(\Omega), \psi_R]$, $[\tr(\Omega), \psi_L]$ and $[\tr(\Omega), \psio_L]$ all vanish, while $[\tr(\Omega), \psio_L]$ is proportional to $\eps$ or $\eps^2\psi_R$ and can be neglected. Hence in all subsequent calculations, we are allowed to replace $\Omega_l^{k} (\psio_k \psi^{l})_L$ by the trace $\tr(\Omega)$ and remember that it commutes with any expression up to negligible terms (in \cite{Per2012} we could even assume that $\tr(\Omega)=0$ by working with a Riemannian connection; this can unfortunately not be done here because of our particular choice of equivariant connection). \\

Then, another Duhamel expansion in the exponentials $\exp(-t_i \nabla^2)$ in the JLO formula, keeping in mind that $\kappa^{2}=0$, yields   
\[\int_{\Delta_{r+1}} \Tr_s\big(e^{-t_{i+1}\nabla^2} [\nabla,\sigma_{i+1}] \ldots  e^{-t_{r+1}\nabla^2} \sigma_0 \,e^{-t_0\nabla^2} [\nabla,\sigma_1] \ldots e^{-t_i\nabla^2}\big) dt \]
as a sum (up to signs) of the following terms :
\begin{gather*}
\int_{\Delta_{r+1}} \Tr_s\big(e^{-t_{i+1}\DD^2} [\nabla,\sigma_{i+1}] \ldots  e^{-t_{r+1}\DD^2} \sigma_0 \,e^{-t_0\DD^2} [\nabla,\sigma_1] \ldots e^{-t_i\DD^2}\big) dt \\
\int_{\Delta_{r+2}} \Tr_s \big(e^{-t_{i+1}\DD^2} \kappa \delta \DD e^{-t_{i+2}\DD^2} [\nabla,\sigma_{i+1}] \ldots  e^{-t_{r+2}\DD^2} \sigma_0 \,e^{-t_0\nabla^2} [\nabla,\sigma_1] \ldots e^{-t_i\DD^2}\big) dt \\
\int_{\Delta_{r+2}} \Tr_s \big(e^{-t_{i+1}\DD^2} [\nabla,\sigma_{i+1}] e^{-t_{i+2}\DD^2} \kappa \delta \DD e^{-t_{i+3}\DD^2} \ldots  e^{-t_{r+2}\DD^2} \sigma_0 \,e^{-t_0\nabla^2} [\nabla,\sigma_1] \ldots e^{-t_i\DD^2}\big) dt \\
\quad \vdots 
\end{gather*}
Notice that the dimension of the simplex on which we integrate is increased by one. These terms also can be rewritten with the modular action $\sigma_{-\DD^{2}}$ : 
\begin{gather*}
\int_{\Delta_{r+1}} \Tr_s\big(\sigma_{-\DD^{2}}^{t_{i+1}}([\nabla,\sigma_{i+1}])\sigma_{-\DD^{2}}^{t_{i+1} + t_{i+2}}([\nabla,\sigma_{i+2}]) \ldots \exp(-\DD^{2})) dt \\ 
\int_{\Delta_{r+2}} \Tr_s \big(\sigma_{-\DD^{2}}^{t_{i+1}}(\kappa \delta \DD) \sigma_{-\DD^{2}}^{t_{i+1} + t_{i+2}}([\nabla,\sigma_{i+1}]) \ldots \exp(-\DD^{2}) \big) dt \\
\quad \vdots 
\end{gather*}
Then for $X = \kappa \delta \DD$, or $\sigma_0$, or a commutator $[\nabla, \sigma_i]$, we have 
\[ \sigma_{-\DD^{2}}^{t}(X) = X + \sum_{k \geq 1} \dfrac{(-t)^{k}}{k!} \ad_{-\DD^{2}}^{k}(X) \]
an so on. We know that $-\DD^{2} \simeq \Delta + \Omega_l^{k} (p_{kL} \partial_{p_l} + (\psio_k \psi^{l})_L)$, and the curvature $\Omega_l^{k} (\psio_k \psi^{l})_L$ can be replaced by the trace $\tr(\Omega)$ which commute with any other expression up to negligible terms. Hence
\begin{gather*}
-[\DD^{2}, X] \simeq \big[\Delta + \Omega_l^{k}p_{kL}\partial_{p_l}, X \big] \simeq \partial_{p_i} X \big(\i \eps \partial_{x^{i}} +  \Omega_i^{k}p_{kL}\big) \\
\big[\DD^2,[\DD^{2}, X]\big] \simeq \partial_{p_i} \partial_{p_j} X \big(\i \eps \partial_{x^{i}} +  \Omega_i^{k}p_{kL}\big)\big(\i \eps \partial_{x^{j}} +  \Omega_j^{l}p_{l L}\big) + \Omega_i^j \partial_{p_i} X \big(\i \eps \partial_{x^{j}} +  \Omega_j^{l}p_{lL}\big) 
\end{gather*}
and continuing the process by induction, we finally get
\[ \sigma_{-\DD^{2}}^{t}(X) \simeq X + \sum_{k \geq 1} \dfrac{t^{k}}{k!} \sum_{\vert \alpha \vert = 1}^{k} P_{\alpha}(X)(\i \eps \partial_{x} + p_{L} \cdot \Omega)^{\alpha} \]
where $P_{\alpha}(X)$ is a linear combination of the $p$-partial derivatives of $X$. The operators $(i \eps \partial_{x} + p_{L} \cdot \Omega)^{\alpha}$ commute with every factor in the graded trace when $\vert \alpha \vert \geq 1$, because $x$-derivatives are dropped. Then, they may be moved to the right in front of $\exp(-\DD^{2})$. Using Example \ref{todd matrix}, we find that these quantities does not contribute, in other words,  
\[ \sigma_{-\DD^{2}}^{t}(X) \simeq X \] 
As a consequence, we can drop the action of $\sigma_{-\DD^{2}}$ in the above calculations and obtain the exact formula       
\begin{multline*}
\int_{\Delta_{r+1}} \Tr_s\big(e^{-t_{i+1}\nabla^2} [\nabla,\sigma_{i+1}] \ldots  e^{-t_{r+1}\nabla^2} \sigma_0 \,e^{-t_0\nabla^2} [\nabla,\sigma_1] \ldots e^{-t_i\nabla^2}\big) dt \\
= \dfrac{1}{(r+1)!} \Tr_s \big([\nabla,\sigma_{i+1}] \ldots [\nabla,\sigma_{r}] \sigma_0  [\nabla,\sigma_1] \ldots [\nabla, \sigma_{i}] \exp(-\DD^{2}) \big) \\
\quad - \dfrac{1}{(r+2)!} \big[ \Tr_s \big(\kappa \delta \DD [\nabla,\sigma_{i+1}] \ldots [\nabla,\sigma_{r}] \sigma_0  [\nabla,\sigma_1] \ldots [\nabla, \sigma_{i}] \exp(-\DD^{2}) \big) \\
\shoveright{\quad  + \Tr_s \big( [\nabla,\sigma_{i+1}] \kappa \delta \DD \ldots [\nabla,\sigma_{r}] \sigma_0  [\nabla,\sigma_1] \ldots [\nabla, \sigma_{i}] \exp(-\DD^2)\big)  + \ldots} \\
\quad + \Tr_s \big( [\nabla,\sigma_{i+1}]  \ldots [\nabla,\sigma_{r}] \sigma_0  [\nabla,\sigma_1] \ldots [\nabla, \sigma_{i}] \kappa \delta \DD \exp(-\DD^2)\big) \big]
\end{multline*}
the factors $\frac{1}{(r+1)!}$ and $\frac{1}{(r+2)!}$ come from the volume of the standard simplex. We are only interested in the terms proportional to $\kappa$. Then, knowing that $-\DD^{2}$ can be replaced by $\Delta + p_L \cdot \Omega\cdot\partial_p + \tr(\Omega)$, using Example \ref{todd matrix} together with the identity $\Td(\Omega) \cdot \exp\tr(\Omega) = \Td(-\Omega)$, and having in mind that the symbols $\sigma_i$ are constant in the direction $EG$, we get 
\begin{multline*}
\dfrac{d}{d \kappa}\int_{\Delta_{r+1}} \Tr_s \big(e^{-t_{i+1}\nabla^2} [\nabla,\sigma_{i+1}] \ldots  e^{-t_{r+1}\nabla^2} \sigma_0 \,e^{-t_0\nabla^2} [\nabla,\sigma_1] \ldots e^{-t_i\nabla^2} \big) dt \\
= \dfrac{1}{(r+1)!}  \barint \big(\langle \langle \delta \sigma_{i+1} [D,\sigma_{i+2}] \ldots [D,\sigma_{r}] \sigma_0  [D,\sigma_1] \ldots [D, \sigma_{i}] \Td(-\Omega) \rangle \rangle [n] \\
\shoveright{\quad - \langle \langle [D, \sigma_{i+1}] \delta\sigma_{i+2} \ldots [D,\sigma_{r}] \sigma_0  [D,\sigma_1] \ldots [D, \sigma_{i}] \Td(-\Omega) \rangle \rangle [n] + \ldots} \\
\shoveright{+ (-1)^{r-1} \langle \langle [D,\sigma_{i+1}] \ldots [D,\sigma_{r}] \sigma_0  [D,\sigma_1] \ldots [D, \sigma_{i-1}] \delta\sigma_{i} \Td(-\Omega) \rangle \rangle [n] \big) }  \\
\quad - \dfrac{1}{(r+2)!} \barint \big(\langle \langle \delta \DD [D,\sigma_{i+1}] \ldots [D,\sigma_{r}] \sigma_0  [D,\sigma_1] \ldots [D, \sigma_{i}] \Td(-\Omega) \rangle \rangle [n] \\
\shoveright{\quad  - \langle \langle [D,\sigma_{i+1}] \delta \DD \ldots [D,\sigma_{r}] \sigma_0  [D,\sigma_1] \ldots [D, \sigma_{i}] \Td(-\Omega) \rangle \rangle [n]  + \ldots} \\
\quad + (-1)^{r}\langle \langle [D,\sigma_{i+1}]  \ldots [D,\sigma_{r}] \sigma_0  [D,\sigma_1] \ldots [D, \sigma_{i}] \delta \DD \, \Td(-\Omega) \rangle \rangle [n] \big)
\end{multline*}

The brackets only select the operators proportional to $(\psio_1 \ldots \psio_n)_R$, which are of order $\leq -(v+2h)$. However, in the quantities proportional to $\frac{1}{(r+1)!}$, these operators gain an extra factor $\delta \sigma$, which is of order $-1$. So, this part is killed by the Wodzicki residue. Moreover, in a local chart $(x,p)$ over $T^*M$ the log-polyhomogeneous symbol $\ln \tilde{q}$ splits as follows :
\[ \ln \tilde{q} = \ln \vert p \vert' + q_0  \]
where $\ln \vert p\vert'$ is constant along $EG$, whereas $q_0$ is a (non-constant) function on $EG$ taking values in the \emph{classical} Heisenberg symbols of order $\leq 0$. Then, 
\begin{align*}
[D, \sigma]& \simeq \i\eps \psi^{i}_R \partial_{x^i} \sigma + \psio_{iR} \partial_{p_i} \sigma  \\
- \delta \DD & \simeq   \big( (d_H q_0)_L + \i\eps \psi^{i}_R(\partial_{x^i} q_0)_L + \psio_{iR}(\partial_{p_i} q_0)_L \big)  \\
& \qquad \qquad + \left(\sum_{i=1}^v \psio_{iR} \left(\dfrac{p_i^{3}}{\vert p \vert'^{4}} \right)_L + \sum_{i=v+1}^n \psio_{iR} \left( \dfrac{p_i}{\vert p \vert'^{4}} \right)_L \right) 
\end{align*}

From now on, let $d=d_H + d_{T^*M}$ denote the total differential over $EG\times T^*M$. The computation of the Wodzicki residue now exactly follows the lines of the proof of Theorems 6.5 and 6.8 of \cite{Per2012}. This basically amounts to make the identifications $\eps\psi_R^{i} \leftrightarrow dx^{i}$ and $\psio_{iR} \leftrightarrow dp_i$, multiply the brackets $\langle\langle \dots \rangle \rangle[n]$ in the formulas above by the standard symplectic volume form $\omega^n / n!$, and compare with the normalization condition $ \langle (\psio_1 \psi^1 \dots \psio_n \psi^n)_R \rangle = (-1)^n $. Keeping track carefully of the numerous signs and factors $\i$, in the present case we get
\begin{multline*}
\dfrac{d}{d \kappa}\int_{\Delta_{r+1}} \Tr_s \big(e^{-t_{i+1}\nabla^2} [\nabla,\sigma_{i+1}] \ldots  e^{-t_{r+1}\nabla^2} \sigma_0 \,e^{-t_0\nabla^2} [\nabla,\sigma_1] \ldots e^{-t_i\nabla^2} \big) dt \\
= \dfrac{1}{(2\pi\i)^n} \dfrac{1}{(r+2)!} \int_{S_H^*M} \iota_L \cdot \Big( \big( \alpha d\rho_{i+1} \ldots d\rho_{r} \rho_0  d\rho_1 \ldots d\rho_{i} - d\rho_{i+1} \alpha d\rho_{r} \rho_0  d\rho_1 \ldots d\rho_{i} + \ldots \\
+ (-1)^{r} d\rho_{i+1}  \ldots d\rho_{r} \rho_0  d\rho_1 \ldots d\rho_{i} \alpha \big) \wedge \Td(R)  \Big)_{\text{vol}} 
\end{multline*}
where $R_l^{k}= d_H\Gamma_{il}^{k} dx^i + \frac{1}{2}R_{lij}^{k} dx^{i}\wedge  dx^{j}$ is the equivariant curvature two-form of the connection, $\rho_i=\rho_*(\hat{a}_i)$, the subscript $_{\text{vol}}$ denotes the component of a differential form of maximal degree with respect to the cotangent bundle, $L$ is the generator of the dilation flow $ F_t(p_1, \dots, p_n) = (e^t p_1, \dots, e^t p_v, e^{2t} p_{v+1}, \dots, e^{2t} p_{n}) $  given in local coordinates by 
\[ L = \left(\sum_{i=1}^v p_i \partial_{p_i} + \sum_{i=v+1}^n 2p_i \partial_{p_i} \right) \]
and $\alpha$ is the one-form on $EG\times T^*M$
\[ \alpha =  dq_0 + \left(\sum_{i=1}^v dp_i \dfrac{p_i^{3}}{\vert p \vert'^{4}} + \sum_{i=v+1}^n dp_i \dfrac{p_i}{\vert p \vert'^{4}}\right) \ .\]
In particular, notice that 
\[ \iota_L \cdot \left(\sum_{i=1}^v dp_{i}\dfrac{p_i^{3}}{\vert p \vert'^{4}} + \sum_{i=v+1}^n dp_{i}\dfrac{p_i}{\vert p \vert'^{4}} \right) = 1 \]
Moreover, the action of the interior product $\iota_L$ on any other one-form $dq_0$ or $d\rho_i$ will not contribute to the Wodzicki residue. Indeed, we already extracted symbols of order $-(v+2h)$, which allows to work at the principal symbol level only. Moreover, these one-forms are (locally) Heisenberg homogeneous functions of order $0$. Thus they are constant along the dilation flow on $T^*M$, so that $\iota_L dq_0 \simeq 0 \simeq \iota_L d\rho_i$. Finally we are left with 
\begin{multline*}
\dfrac{d}{d \kappa}\int_{\Delta_{r+1}} \Tr_s \big(e^{-t_{i+1}\nabla^2} [\nabla,\sigma_{i+1}] \ldots  e^{-t_{r+1}\nabla^2} \sigma_0 \,e^{-t_0\nabla^2} [\nabla,\sigma_1] \ldots e^{-t_i\nabla^2} \big) dt \\
= \dfrac{1}{(2\pi\i)^n}\dfrac{1}{(r+1)!} \int_{S_H^*M} d\rho_{i+1} \ldots d\rho_{r} \rho_0  d\rho_1 \ldots d\rho_{i} \wedge \Td(R)
\end{multline*}

Analogous manipulations on the second sum of the cochain $\chi^{\Tr_s}(\sigma_{*}, \nabla)$ give the final answer. $\hfill{\square}$ 
\end{pr}

Propositions \ref{pcocycle}, \ref{pDeRham} and \ref{pAffine} give the following sequence of cohomologous cocycles in the complex $\Hom(\Omh\Th\A,\C)$ :
\[ \partial([\tau])  \equiv \lambda'_1 \circ \chi^{\Tr_s}(\sigma_*,\DD,\ln \tilde{q}_L) \equiv \lambda'_{\Td(R)} \circ \chi(\rho_*,d) \equiv \Phi(\Td(TM\otimes\C)) \ . \]
We thus get the main result of this paper :

\begin{thm}\label{ttodd}
Let $M$ be a foliated manifold, $G\subset\Diff(M)$ a discrete group of diffeomorphisms mapping leaves to leaves, $0\to \Psi^{-1}_{H,c}(M)\rtimes G \to \Psi^0_{H,c}(M)\rtimes G \to C^{\infty}_c (S^*_HM) \rtimes G \to 0$ the equivariant Heisenberg pseudodifferential extension. Then the image of the canonical trace localized at unit $[\tau]\in \HP^0 (\Psi^{-1}_H(M)\rtimes G)$ under the excision map is
\begin{equation}
\partial([\tau]) = \Phi(\Td(TM\otimes\C)) 
\end{equation}
where $\Phi: H^{\mathrm{ev}}(EG\times_G S^*_HM) \to HP^1(C^{\infty}_c(S^*_HM)\rtimes G)$ is Connes' characteristic map from equivariant cohomology to cyclic cohomology, and $\Td(TM\otimes\C) $ is the equivariant Todd class of the complexified tangent bundle of $M$.
\end{thm}

\section{The transverse index theorem of Connes and Moscovici} \label{CM signature}

Let $\A$ be an associative algebra and $(H,F)$ a (trivially graded) $p$-summable Fredholm module. Hence, $\A$ is represented by bounded operators on a separable Hilbert space $H$, and $F$ is a bounded self-adjoint operator on $H$ such that the operators $a(F^2-1)$, $(F^2-1)a$ and $[F,a]$ are in the Schatten class $\ell^p(H)$ for all $a\in A$. In addition, we suppose given an extension of ``abstract pseudodifferential operators''
\begin{equation}
0 \to \Psi^{-1} \to \Psi^0 \to \Psi^0/\Psi^{-1} \to 0 \label{ext}
\end{equation}
where 
\begin{itemize}
\item $\Psi^0$ is an algebra of bounded operators on $H$ containing the representation of $\A$,
\item $\Psi^{-1}$ is a two-sided ideal consisting of $p$-summable operators on $H$, 
\item $F$ is a multiplier of $\Psi^0$ and $[F,\Psi^0]\subset \Psi^{-1}$.
\end{itemize}
Let $P=\frac{1}{2}(1+F)$. Then $[P,a]\in \Psi^{-1}$ and $aP^2\equiv aP\mod \Psi^{-1}$ for all $a\in A$. The linear map
\begin{equation}
\rho_F\ :\ \A \to \Psi^0/\Psi^{-1}\ ,\qquad \rho_F(a) \equiv aP\mod \Psi^{-1}\ ,
\end{equation}
is an algebra homomorphism since $a_1Pa_2P \equiv a_1a_2 P\mod \Psi^{-1}$ for all $a_1,a_2\in \A$.
\begin{lem}
The Chern-Connes character of the Fredholm module $(H,F)$ is given by the odd cyclic cohomology class over $\A$
\begin{equation}
\ch(H,F) =\rho^*_F \circ \partial([\Tr]) 
\end{equation}
where $[\Tr]\in HP^0(\Psi^{-1})$ is the class of the operator trace, $\partial:HP^0(\Psi^{-1}) \to HP^1(\Psi^0/\Psi^{-1})$ is the excision map associated to extension (\ref{ext}), and $\rho_F^*: HP^1(\Psi^0/\Psi^{-1})\to HP^1(\A)$ is induced by the homomorphism $\rho_F$. 
\end{lem}
\begin{pr}
Consider the algebra $\E = \{ (Q,a) \in \Psi^0\oplus \A \ |\ Q\equiv aP\mod \Psi^{-1}\}$. The homomorphism $\E \to\A$, $(Q,a)\mapsto a$ yields an extension
\[
0\to \Psi^{-1} \to \E \to \A \to 0\ .
\]
By definition (\cite{ConIHES}), the Chern-Connes character $\ch(H,F)\in HP^1(\A)$ is the image of the operator trace under the excision map associated to this extension. On the other hand, the homomorphism $\E\to\Psi^0$, $(Q,a)\mapsto Q$ yields a commutative diagram of extensions
$$
\xymatrix{ 0 \ar[r] & \Psi^{-1} \ar@{=}[d] \ar[r] & \E \ar[d] \ar[r] & \A \ar[d]^{\rho_F} \ar[r] & 0 \\
0 \ar[r] & \Psi^{-1} \ar[r] & \Psi^0 \ar[r] & \Psi^0/\Psi^{-1} \ar[r] & 0 }
$$
The conclusion then follows from the naturality of excision. $\hfill{\square}$ 
\end{pr}
We apply this to the hypoelliptic operators constructed by Connes and Moscovici in \cite{CM1995}. Let $M$ be an oriented foliated manifold and $G\subset\Diff(M)$ a discrete group of orientation-preserving diffeomorphisms mapping leaves to leaves. We make the hypothesis that $G$ has \emph{freely}, so that only the identity element has fixed points on $M$. We denote by $V\subset TM$ the subbundle tangent to the leaves and by $N=TM/V$ the normal bundle; both are equivariant $G$-bundles by construction. Assume that $V$ and $N$ are provided with $G$-invariant euclidean structures, called $G$-invariant \emph{triangular structures} in \cite{CM1995}. Then the hermitean vector bundle
\begin{equation}
E = \Lambda^\bullet(V^*\otimes\C)\otimes \Lambda^\bullet(N^*\otimes\C)
\end{equation}
is $G$-equivariant, and the euclidean structures on $V,N$ determine a $G$-invariant volume form on $M$ via the canonical isomorphism of top-degree forms $\Lambda^{\max}V\otimes\Lambda^{\max}N \cong \Lambda^{\max}M$. Let $H = L^2(M,E)$ be the Hilbert space of square-integrable sections of $E$ with respect to the hermitean structure and volume form. The crossed-product algebra 
\begin{equation}
\A = C^{\infty}_c(M)\rtimes G
\end{equation} 
is represented by bounded operators on $H$ as follows: a function $f\in C^{\infty}_c(M)$ acts on the sections of $E$ by pointwise multiplication, while $g\in G$ is represented by the unitary operator coming from the action of $G$ on the manifold $M$ and the vector bundles $V,N$. Denote by $d_V: C^{\infty}(M,E)\to C^{\infty}(M,E)$ the leafwise de Rham differential. Choose an isomorphism of $N$ with a vector subbundle of $TM$ transverse to $V$, and denote by $d_N$ the corresponding transverse de Rham differential. Then Connes and Moscovici consider the \emph{hypoelliptic signature operator} acting on $C^{\infty}(M,E)$
\begin{equation}
Q = \pm(d_Vd_V^*-d_V^*d_V) + (d_N+d_N^*)\ ,
\end{equation}
where the sign $+1$ is taken on $\Lambda^{\mathrm{ev}}N^*$ and $-1$ on $\Lambda^{\mathrm{odd}}N^*$. This is a formally self-adjoint, hypoelliptic differential operator of order two. $Q$ is not quite invariant under the action of $G$ because the isomorphism $TM \cong V\oplus N$ requires a choice. However, in the Heisenberg pseudodifferential calculus associated to the foliation on $M$, the operator $Q$ is Heisenberg-elliptic and its leading symbol is exactly $G$-invariant. From this one can build a properly supported Heisenberg pseudodifferential operator
\begin{equation}
F = \frac{Q}{|Q|}
\end{equation}
which is defined only up to addition of a smoothing operator. Again the Heisenberg leading symbol of $F$ is $G$-invariant. Now we turn to the geometric example of \cite{CM1995}, where $M$ is the bundle of Riemannian metrics over a smooth $G$-manifold $W$. Here the foliation on $M$ corresponds to the fibration $M\to W$, and has a tautological triangular structure. The action of $G$ by diffeomorphisms on $W$ canonically lifts to an action on $M$ mapping leaves to leaves, and preserving the triangular structure. In this situation, the results of \cite{CM1995} show that $F$ is a bounded operator on $H$, $a(F^2-1)$ and $(F^2-1)a$ are smoothing for all $a\in\A$, and the pair $(H,F)$ defines a $p$-summable Fredholm module over the algebra $\A$ for any $p> \dim V +2\dim N$. Its Chern-Connes character may thus be computed by means of the above lemma. We let $\Psi_{H,c}(M,E)$ be the algebra of compactly supported Heisenberg pseudodifferential operators acting on the smooth sections of $E$. Let $\pi: S^*_HM \to M$ be the projection from the Heisenberg cosphere bundle. The pullback $\pi^*E$ is naturally a $G$-equivariant vector bundle over $S^*_HM$, and $C^{\infty}_c(S^*_HM,\End(\pi^*E))$ is the $G$-algebra of Heisenberg leading symbols. The representation of the crossed-product $\Psi_{H,c}(M,E)\rtimes G$ on the Hilbert space $H$ leads to subalgebras of bounded operators
\[
\Psi^{0} = \Im(\Psi_{H,c}^{0}(M,E)\rtimes G)\ ,\qquad \Psi^{-1} = \Im(\Psi^{-1}_{H,c}(M,E)\rtimes G)\ .
\]
These representations are not faithful. Observe that any operator in $\Psi^0$ has a smooth Schwartz kernel on $M\times M$ except maybe at the points $(x,x\cdot g)$, where $(x,g)\in M\times G$. Since by hypothesis $G$ acts freely on $M$, the singular set of the Schwartz kernel is contained in a disjoint union of submanifolds diffeomorphic to $M$ in $M\times M$, indexed by the elements $g\in G$. Therefore the leading symbol map $\Psi^0_{H,c}(M,E)\rtimes G\to C^{\infty}_c(S^*_HM,\End(\pi^*E))\rtimes G$ factors through $\Psi^0$, and conversely a leading symbol completely determines the class of an operator in $\Psi^0/\Psi^{-1}$. One thus has a canonical isomorphism of algebras
\[
\Psi^0/\Psi^{-1} \cong C^{\infty}_c(S^*_HM,\End(\pi^*E))\rtimes G\ .
\]
Under this identification the homomorphism $\rho_F: \A \to \Psi^0/\Psi^{-1}$ is given by
\[
\rho_F(fU_g) = \pi^*(f)e U^E_g\qquad \forall\ f\in C^{\infty}_c(M)\ ,\ g\in G
\]
where $e\in  C^{\infty}(S^*_HM,\End(\pi^*E))$ is the leading symbol of the operator $P=\frac{1}{2}(1+F)$, and $U^E_g$ is the represntation of $g$ as a linear operator on the space of sections of $\pi^*E$. Since $P^2\equiv P$ and $PU^E_g\equiv U^E_gP$ modulo operators of order $-1$, one has $e^2=e$ and $eU^E_g=U^E_ge$ for all $g\in G$. Hence $e$ is a $G$-invariant idempotent section of the bundle $\End(\pi^*E)$. Its range is the $G$-equivariant subbundle $E_+$ of $\pi^*E$ consisting in the positive eigenvectors for the leading symbol of $F$. By Chern-Weil theory, the equivariant Chern character $\ch(E_+)$ is represented by a closed $G$-invariant differential form on the homotopy quotient $EG\times_G S^*_HM$. Taking its product with the equivariant Todd class of the complexified tangent bundle yields a class
\begin{equation}
L'(M) = \Td(TM\otimes\C) \cup \ch(E_+) \in H^{\mathrm{ev}}(EG\times_G S^*_HM)\ .
\end{equation}

\begin{thm}
Let $G$ be a discrete group of orientation-preseving diffeomorphisms on a smooth oriented manifold $W$. Let $M$ be the bundle of Riemannian metrics over $W$ and $\A = C^{\infty}_c(M)\rtimes G$. If $G$ has no fixed points, then the Chern-Connes character of the Fredholm module $(H,F)$ associated to the hypoelliptic signature operator of Connes and Moscovici is   
\begin{equation}
\ch(H,F) = \pi_*\circ\Phi(L'(M))\ \in HP^1(\A)\ ,
\end{equation}
where $\Phi: H^{\mathrm{ev}}(EG\times_G S^*_HM) \to HP^1(C^{\infty}_c(S^*_HM)\rtimes G)$ is Connes' characteristic map from equivariant cohomology to cyclic cohomology, and $\pi_*: HP^1(C^{\infty}_c(S^*_HM)\rtimes G) \to HP^1(\A)$ is the map induced by the projection $\pi:S^*_HM\to M$.
\end{thm}
\begin{pr}
One has to compare the two extensions
$$
\xymatrix{ 0 \ar[r] & \Psi^{-1}_{H,c}(M,E)\rtimes G \ar[d] \ar[r] & \Psi^{0}_{H,c}(M,E)\rtimes G \ar[d] \ar[r] & C^{\infty}_c(S^*_HM,\End(\pi^*E))\rtimes G \ar@{=}[d] \ar[r] & 0 \\
0 \ar[r] & \Psi^{-1} \ar[r] & \Psi^0 \ar[r] & \Psi^0/\Psi^{-1} \ar[r] & 0 }
$$
where the vertical arrows are the representations as bounded operators in the Hilbert space $H$. We consider two different cyclic cohomology classes on the ideals. The first one is the operator trace $[\Tr]\in HP^0(\Psi^{-1})$, and the second is the trace \emph{localized at unit} $[\tau]\in HP^0(\Psi^{-1}_{H,c}(M,E)\rtimes G)$. Of course $[\tau]$ is not the pullback of $[\Tr]$ under the representation. We use a zeta-function renormalization in order to compute the image $\partial([\Tr])\in HP^1(\Psi^0/\Psi^{-1})$ of the operator trace under the excision map of the bottom extension, as in Section \ref{sindex}. Then, \emph{since $G$ has no fixed points}, only the part of the operator trace which is localized at units contributes to the residues. This means that one has the equality
\[
\partial ([\Tr]) = \partial([\tau])
\]
in $HP^1(C^{\infty}_c(S^*_HM,\End(\pi^*E))\rtimes G)$. A choice of local trivializations of the vector bundle $E$ and a partition of unity allows to identify $C^{\infty}_c(S^*_HM,\End(\pi^*E))\rtimes G$ with a subalgebra of the algebra of matrices $M_{\infty}(C^{\infty}_c(S^*_HM)\rtimes G)$. Under this identification Theorem \ref{ttodd} implies the equality
\[
\partial([\Tr]) = \tr\#\Phi(\pi^*\Td(TM\otimes\C))
\]
where $\tr$ denotes the trace on $M_{\infty}(\C)$ and $\#$ is the cup-product of cyclic cocycles (\cite{ConIHES}). Finally the homomorphism $\rho_F$ is multiplication by the $G$-invariant idempotent $e\in \cinf(S^*_HM,\End(\pi^*E))\subset M_{\infty}(\cinf(S^*_HM))$, so the composition $\rho_F^*\circ\partial([\Tr])$ is the above class twisted by the Chern character $\ch(E_+)$. $\hfill{\square}$ 
\end{pr}


\end{document}